\documentclass[11pt,reqno,a4paper]{amsart}

\usepackage{amsmath}
\usepackage{enumerate}
\usepackage{amssymb}
\usepackage{amscd}
\usepackage{amsthm}
\usepackage{amsfonts}
\usepackage{graphicx}
\usepackage{hyperref}
\usepackage[matrix, arrow, curve]{xy}
\usepackage{etoolbox}
\usepackage{comment}

\usepackage{fouriernc}
\usepackage[T1]{fontenc}

\numberwithin{equation}{section}

\patchcmd{\subsection}{-.5em}{.5em}{}{}
\patchcmd{\subsubsection}{-.5em}{.5em}{}{}

\newcommand{\bA}{\mathbb{A}}

\newcommand{\bC}{\mathbb{C}}
\newcommand{\bD}{\mathbb{D}}

\newcommand{\bH}{\mathbb{H}}

\newcommand{\bR}{\mathbb{R}}

\newcommand{\R}{\mathbb{R}}
\newcommand{\C}{\mathbb{C}}

\newcommand{\ra}{\rightarrow}

\newcommand{\qand}{\quad \textrm{and} \quad}

\newcommand\subsetsim{\mathrel{%
\ooalign{\raise0.2ex\hbox{$\subset$}\cr\hidewidth\raise-0.8ex\hbox{\scalebox{0.9}{$\sim$}}\hidewidth\cr}}}

\DeclareMathOperator{\supp}{supp}
\DeclareMathOperator{\trace}{tr}

\DeclareMathOperator{\arcosh}{arcosh}
\DeclareMathOperator{\artanh}{artanh}

\newcommand{\Q}{\mathbb Q}
\newcommand{\Z}{\mathbb Z}

\renewcommand{\epsilon}{\varepsilon}

\theoremstyle{theorem}
\newtheorem{theorem}{Theorem}[section]
\newtheorem{corollary}[theorem]{Corollary}
\newtheorem{proposition}[theorem]{Proposition}
\newtheorem{lemma}[theorem]{Lemma}

\theoremstyle{definition}
\newtheorem{definition}[theorem]{Definition}
\newtheorem{convention}[theorem]{Convention}
\newtheorem{remark}[theorem]{Remark}
\newtheorem{caveat}[theorem]{Caveat}
\newtheorem*{example}{Example}

\renewcommand{\phi}{\varphi}

\begin{document}

\title[Aperiodic order and spherical diffraction]{Aperiodic order and spherical diffraction, II:\\ Translation bounded measures on homogeneous spaces}

%  Author I information
\author{Michael Bj\"orklund}
\address{Department of Mathematics, Chalmers, Gothenburg, Sweden}
\email{micbjo@chalmers.se}
\thanks{}

%    Author II information
\author{Tobias Hartnick}
\address{Institut für Algebra und Geometrie, KIT, Karlsruhe, Germany}
\curraddr{}
\email{tobias.hartnick@kit.edu}
\thanks{}

\author{Felix Pogorzelski}
\address{Mathematisches Institut, Universität Leipzig, Leipzig, Germany}
\curraddr{}
\email{felix.pogorzelski@math.uni-leipzig.de}
\thanks{}


\subjclass[2010]{Primary: ; Secondary: }

\date{}

\dedicatory{}

\maketitle

\begin{abstract} We study the auto-correlation measures of invariant random point processes in the hyperbolic plane which arise from various classes of aperiodic Delone sets. More generally, we study auto-correlation measures for large classes of Delone sets in (and even translation bounded measures on) arbitrary locally compact homogeneous metric spaces. We then specialize to the case of weighted model sets, in which we are able to derive more concrete formulas for the auto-correlation. In the case of Riemannian symmetric spaces we also explain how the auto-correlation of a weighted model set in a Riemannian symmetric space can be identified with a (typically non-tempered) positive-definite distribution on $\R^n$. This paves the way for a diffraction theory for such model sets, which will be discussed in the sequel to the present article. 
\end{abstract}

\section{Introduction}

\subsection{General themes of this article}
The study of aperiodic Delone sets in $\R^n$ and more general locally compact abelian groups is a classical topic in harmonic analysis (see \cite{BaakeGrimm} for an extensive reference list). A particular interesting class of such Delone sets are model sets as introduced by Meyer in his pioneering work \cite{Meyer}. In the first part of this series of articles \cite{BHP1} we have studied model sets in the wider setting of - typicall non-abelian -  locally compact second-countable (lcsc) groups and developed a theory of auto-correlation for such model sets (and more generally, for so-called Delone sets of finite local complexity in locally compact groups).

In this second part we study Delone sets, i.e.\ uniformly discrete and relatively dense subsets, in arbitrary lcsc homogeneous metric spaces. Here a locally compact metric space is called \emph{homogeneous} if its isometry group $G$ acts transitively on $X$. Examples of such spaces exist in abundance; we will consider in particular Euclidean spaces, hyperbolic spaces, Riemannian symmetric spaces, vertex sets of regular trees and Bruhat--Tits buildings and locally compact second countable (lcsc) groups themselves with an invariant metric. 

Any locally compact homogeneous metric space is of the form $X = K\backslash G$ for a compact subgroup $K<G$ and we will show that every Delone set in $X$ is the orbit of a Delone set in $G$ as defined in \cite{BH}. In particular, we can define a \emph{model set} in $X$ as the orbit of a model set in $G$, and these are the main protagonists of the current article. While the case of Euclidean space (seen as homogeneous space under the group of Euclidean motions, \cite{BFG}) and abelian locally compact groups \cite{BaakeGrimm} have been studied before, this seems to be the first systematic investigation of auto-correlation of Delone sets in general lcsc homogeneous metric spaces.

The bulk of this article is devoted to transferring the theory of auto-correlation developed in \cite{BHP1} from model sets in lcsc groups to model sets (and more general translation bounded measures) in arbitrary lcsc homogeneous metric spaces. 

While our results apply in large generality, this introduction will focus on the simple special case of Delone sets in the hyperbolic plane, for which we can state some of our results in a particularly nice form. In particular we are going to explain how a model set in the hyperbolic plane gives rise to an evenly positive-definite (generally non-tempered) distribution on the real line. The complex Fourier transform of this distribution will be the subject of the third paper in this series \cite{BHP3}, where it will be established that it is a pure point Radon measure if the model set is uniform. The natural context of this result is the theory of spherical diffraction alluded to in the title of this series of articles.

\subsection{Tilings of the hyperbolic plane}

\begin{figure}[ht]
	\centering
  \includegraphics[width=\columnwidth]{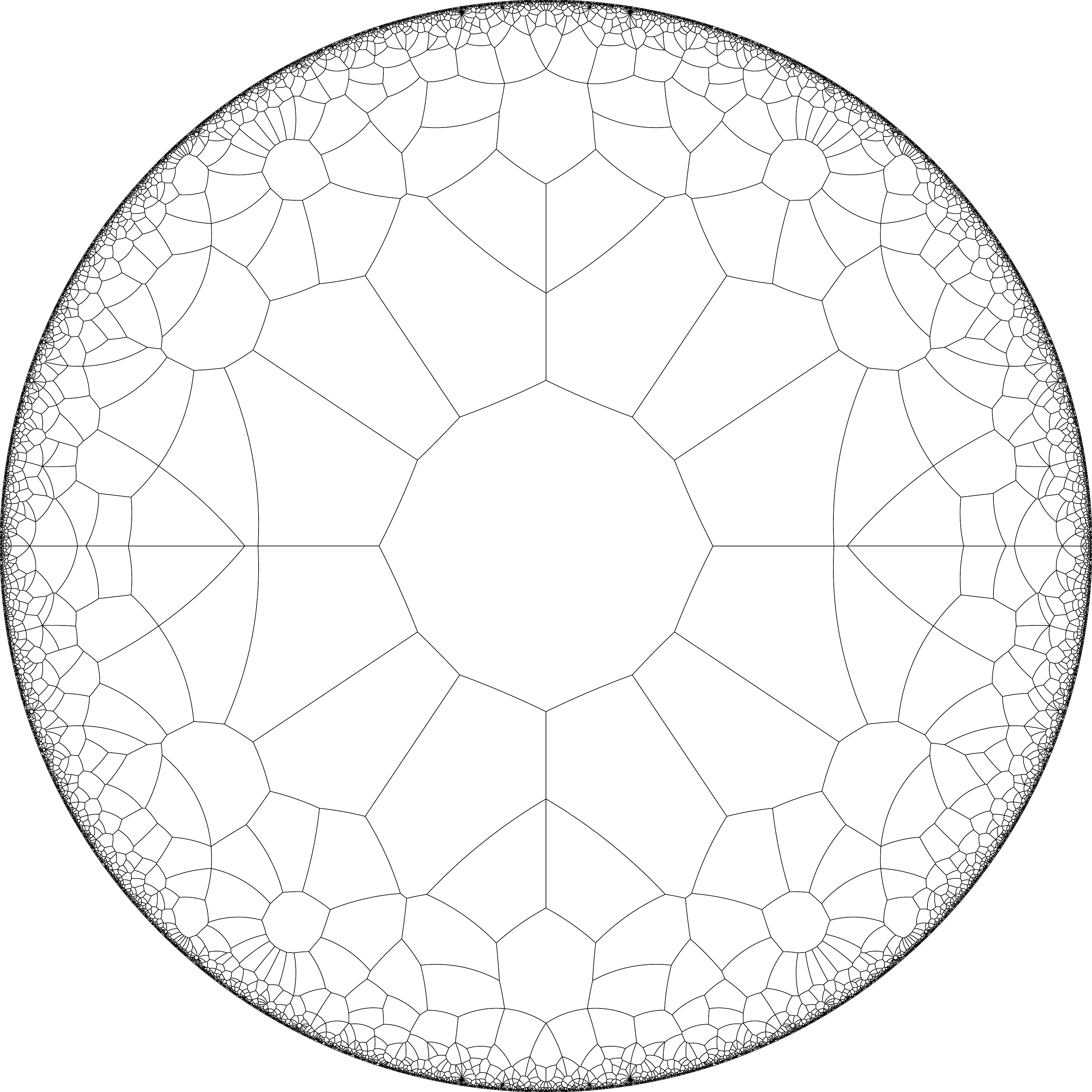}
	\caption{Voronoi tiling of a hyperbolic model set	 (Picture courtesy of Stefan Witzel)}
	\label{Figure1}
\end{figure}

The \emph{Poincar\'e disc model} of the hyperbolic plane is given by the unit disc $\bD \subset \bC$ with the metric
\[
d(z_1, z_2) := 2 \artanh\left(\frac{|z_1-z_2|}{|1-z_1\overline{z_2}|}\right)
%\left(\frac{|z_1-z_2|}{|z_1-\overline{z_2}|}\right)
\]
A subset $\Lambda \subset \bD$ is called a \emph{Delone set} if it is uniformly discrete and relatively dense, i.e.\ if there exist constants $R>r>0$ such that $d(\lambda_1, \lambda_2) \geq r$ for all $\lambda_1, \lambda_2 \in \Lambda$ with $\lambda_1 \neq \lambda_2$ and if for every $z \in \bD$ there exists $\lambda \in \Lambda$ with $d(\lambda, z) \leq R$. It is called \emph{periodic} if the group $\Gamma := \{g \in {\rm Is}(\bD, d) \mid g.\Lambda = \Lambda \}$ acts cocompactly on $\bD$.

If $\Lambda$ is a Delone set in the hyperbolic plane, then the \emph{Voronoi cell} of $\lambda \in \Lambda$ is the convex compact set with piecewise-geodesic boundary given by
\[
V_\lambda = \{z \in \bH^2 \mid \forall \lambda' \in \Lambda \setminus \{\lambda\}:\, d(z, \lambda) \leq d(z, \lambda')\}.
\]
The Voronoi cells $(V_\lambda)_{\lambda \in \Lambda}$ form a tiling of the hyperbolic plane called the \emph{Voronoi tiling} of $\Lambda$. Figure \ref{Figure1}, due to Stefan Witzel, shows a piece of a Voronoi tiling associated with a Delone set in the Poincar\'e disc. The underlying Delone set is not periodic, but nevertheless enjoys a great deal of structure, analogous to ``quasi-crystals'' in the Euclidean plane. In fact, in the terminology introduced below, it is a uniform model set in the Poincar\'e disc.

We remark that the study of such non-periodic tilings in the hyperbolic plane has a long history (see e.g.\ \cite{M97,MM98,BW92,BH13}), but we will see that hyperbolic model sets and their associated tilings have a number of features which are not known to hold in previous examples.

% and the corresponding Voronoi tiling  %If $\Gamma$ is any uniform lattice in ${\rm Is}(\bD, d)$, i.e.\ a discrete cocompact subgroup, then the orbit 
%\[
%\Gamma.o = \{\gamma(o) \mid \gamma \in \Gamma\}
%\]
%defines a periodic Delone set in $\bD$ and a corresponding Voronoi tiling. 
%In this article we are interested in certain non-periodic Delone sets in the hyperbolic plane, whose associated tilings still enjoy a large amount of order and which are analogous to ``quasi-crystals'' in the Euclidean plane. 

%The study of such non-periodic tilings in the hyperbolic plane has a long history with fundamental contributions by Margulis, Mozes and others (REFERENCES). Our construction of non-periodic Delone sets in the hyperbolic sets will be based on the notion of model sets in non-abelian groups as developed in \cite{BHP1}, and we will see that these Delone sets have a number of remarkable properties not shared by previous examples.

\subsection{Weighted model sets and unique ergodicity}\label{IntroWMS} In the sequel it will be convenient for us to work with the \emph{upper half-plane} model of the hyperbolic plane as given by
\[
\bH^2 = \{(x,y) \in \R^2 \mid y >0\} = \{z \in C\mid {\rm Im}(z) > 0\}
\]
with metric
\[
d((x_1, y_1), (x_2, y_2)) = \arcosh\left(1 + \frac{(x_1-x_2)^2+(y_1-y_2)^2}{2y_1y_2}\right).
\]
The group ${\rm SL}_2(\R)$ acts on $\bH^2$ by isometries via
\[
\begin{pmatrix} a&b\\c&d \end{pmatrix}.z := \frac{az+b}{cz+d} \quad \left(\begin{pmatrix} a&b\\c&d \end{pmatrix} \in {\rm SL}_2(\R), \, z\in \C,\, {\rm Im}(z)>0\right),
\]
and every orientation-preserving isometry of $\bH^2$ arises from a matrix in ${\rm SL}_2(\R)$ in this way. Moreover, the action ${\rm SL}_2(\R) \curvearrowright \bH^2$ is transitive. 
%and 
%the stabilizer of $i$ is given by ${\rm SO}_2(\R)$, hence we may identity $\bH^2$ with the homogeneous space ${\rm SL}_2(\R)/{\rm SO}_2(\R)$. (For technical reasons we will actually prefer to work with the left-quotient ${\rm SO}_2(\R)\backslash {\rm SL}_2(\R)$, which is however isomorphic.)
\begin{lemma}[Lifting lemma]\label{IntroLift} For every Delone set $\Lambda \subset \bH^2$ there exists a Delone set $\widetilde{\Lambda}$ in ${\rm SL}_2(\R)$ whose orbit coincides with $\Lambda$, i.e.\
\begin{equation}\label{LambdaOrbit}
\Lambda = \{\widetilde{\lambda}.i \mid \widetilde{\lambda} \in \widetilde{\Lambda}\}.
\end{equation}
\end{lemma}
Not every orbit of a Delone set in ${\rm SL}_2(\R)$ is a Delone set in the hyperbolic plane, but if $\Gamma < {\rm SL}_2(\R)$ is a uniform lattice, then the orbit $\Gamma.i$ of $\Gamma$ in $\bH^2$ defines a periodic Delone set in $\bH$. 

We now explain how to construct non-periodic examples of highly structured Delone set in the hyperbolic plane: We start from a uniform lattice $\Gamma$ in ${\rm SL}_2(\R) \times {\rm SL}_2(\R)$, for example
\[
\Gamma = \{(Z_1(x,y), Z_2(x,y)) \mid x_0, \dots, x_3, y_0, \dots, y_3 \in \Z, \det(Z_1(x,y)) = \det(Z_2(x,y)) = 1\},
\]
where
\begin{tiny}
\[
Z_1(x, y) :=  \left(\begin{matrix} \left(x_0 +  \frac{1+\sqrt 5}{2} \cdot y_0\right) + \left(x_1+  \frac{1+\sqrt 5}{2} \cdot y_1\right)  \sqrt {3+\sqrt 5} & \left(x_2 +  \frac{1+\sqrt 5}{2} \cdot y_2\right)+\left(x_3+  \frac{1+\sqrt 5}{2} \cdot y_3\right) \sqrt {3+\sqrt 5} \\  \frac{1-3\sqrt 5}{2}\left( \left(x_2 +  \frac{1+\sqrt 5}{2} \cdot y_2\right) - \left(x_3+  \frac{1+\sqrt 5}{2} \cdot y_3\right)  \sqrt {3+\sqrt 5}\right) &  \left(x_0 +  \frac{1+\sqrt 5}{2} \cdot y_0\right)- \left(x_1+  \frac{1+\sqrt 5}{2} \cdot y_1\right) \sqrt {3+\sqrt 5}\end{matrix}\right),
\]
\[
Z_2(x,y) :=  \left(\begin{matrix}  \left(x_0 +  \frac{1-\sqrt 5}{2} \cdot y_0\right) +  \left(x_1 +  \frac{1-\sqrt 5}{2} \cdot y_1\right)  \sqrt {3-\sqrt 5} &  \left(x_2 +  \frac{1-\sqrt 5}{2} \cdot y_2\right)+ \left(x_3 +  \frac{1-\sqrt 5}{2} \cdot y_3\right) \sqrt {3-\sqrt 5} \\  \frac{1+3\sqrt 5}{2}\left( \left(x_2 +  \frac{1-\sqrt 5}{2} \cdot y_2\right)-  \left(x_3 +  \frac{1-\sqrt 5}{2} \cdot y_3\right)\sqrt {3-\sqrt 5}\right) &  \left(x_0 +  \frac{1-\sqrt 5}{2} \cdot y_0\right)- \left(x_1 +  \frac{1-\sqrt 5}{2} \cdot y_1\right)\sqrt {3-\sqrt 5}\end{matrix}\right).
\]
\end{tiny}
Moreover, let $W$ be a compact identity neightbourhood in ${\rm SL}_2(\R)$ and denote by $p_1, p_2: {\rm SL}_2(\R) \times {\rm SL}_2(\R) \to {\rm SL}_2(\R)$ the two coordinate projections. We then define a subset of ${\rm SL}_2(\R)$ by 
\[
\widetilde{\Lambda} := p_1(\Gamma \cap ({\rm SL}_2(\R) \times W)).
\]
This is an example of a \emph{uniform model set} in ${\rm SL}_2(\R)$, and such model sets have been systematically studied in \cite{BH, BHP1}.  The following proposition holds for model sets and more generally for Delone set in ${\rm SL}_2(\R)$ of finite local complexity as defined in \cite{BH, BHP1}.
\begin{proposition}\label{IntroLift2} If $\widetilde{\Lambda}$ is a Delone set of finite local complexity in ${\rm SL}_2(\R)$, for example a uniform model set, then its orbit $\Lambda$ as defined by \eqref{LambdaOrbit} is a Delone set in $\bH^2$, and for every $\lambda \in \Lambda$ we have
\[
w(\lambda) := |\{\widetilde{\lambda} \in \widetilde{\Lambda} \mid \widetilde{\lambda}.i  = \lambda\}| < \infty.
\]
\end{proposition}
We refer to the pair $(\Lambda, w)$ arising from this construction as a \emph{weighted uniform model set} in $\bH^2$ and to the function $w:\Lambda \to \mathbb N$ as its \emph{weight function}. The associated \emph{weighted Dirac comb} is the Radon measure $\delta_{(\Lambda, w)}$ on $\bH^2$ given by
\[
\delta_{(\Lambda, w)}(f) = \sum_{\lambda \in \Lambda} w(\lambda) f(\lambda) \quad (f \in C_c(\bH^2)).
\]
The weighted Dirac comb of a weighted model set is an example of a translation bounded measure on $\bH^2$. Such measures have been studied extensively in the setting of abelian groups \cite{TbBook, Hof, BaakeLenz}, and generalizing results from the abelian case we will show:
\begin{proposition}\label{HullCompactIntro} The weak-$*$-closure $\Omega_{(\Lambda, w)} := \overline{{\rm SL}_2(\R).\delta_{(\Lambda, w)}}$ in the space of Radon measures on $\bH^2$ compact.
\end{proposition}
One can show that $\Omega_{(\Lambda, w)}$ consists of those weighted Dirac combs of weighted Delone sets in $\bH^2$ which locally coincide with $(\Lambda, w)$ up to an element of ${\rm SL}_2(\R)$. We refer to $\Omega_{(\Lambda, w)}$ as the \emph{hull} of $(\Lambda, w)$. By construction the group ${\rm SL}_2(\R)$ acts on the hull, and we can extend the results from \cite{BHP1} to show:
\begin{theorem}[Unique ergodicity of regular weighted model sets]\label{UniqueErgodicityIntro} If $\widetilde{\Lambda}$ is a regular uniform model set, then the hull $\Omega_{(\Lambda, w)}$ is minimal and admits a unique ${\rm SL}_2(\R)$-invariant probability measure.
\end{theorem}
We refer the reader to \cite{BHP1} for the precise definition of a \emph{regular} uniform model set. Besides some technical conditions on the window it requires the window to be in general position with respect to the lattice $\Gamma$.
%To the best of our knowledge, these are the first examples of non-periodic uniquely ergodic Delone sets in the hyperbolic plane.
\subsection{Auto-correlation measures and auto-correlation distributions}
From now on we fix a weighted uniform model set $(\Lambda, w)$ in the hyperbolic plane arising from a \emph{regular} uniform model set in ${\rm SL}_2(\R)$.

By Theorem \ref{UniqueErgodicityIntro} there exists a unique ${\rm SL}_2(\R)$-invariant  probability measure on the hull $\Omega_{(\Lambda, w)}$. If we denote this measure by $\nu$, then the pair $(\Omega_{(\Lambda, w)}, \nu)$ is an example of a (weighted) point process in the hyperbolic plane, and for such point processes one can define correlation measures in the usual way. For example, the two-point correlation $\eta^{(2)}$ is the Radon measure on $\bH^2 \times \bH^2$ given by
\[
\eta^{(2)}(f_1 \otimes f_2) = \int_{\Omega_{(\Lambda, w)}} \left(\int_{\bH^2} f_1 \, d\mu\right) \left(\int_{\bH^2} f_2 \, d\mu\right) d\nu(\mu).
\]
Since $\nu$ is ${\rm SL}_2(\R)$-invariant, the two-point correlation descends to a Radon measure $\eta$ on the quotient ${\rm SL}_2(\R)\backslash (\bH \times \bH)$ called the \emph{auto-correlation measure} of $(\Lambda, w)$ (or of $\nu$).

There are several ways to think of this measure. Firstly, if we abbreviate $G := {\rm SL}_2(\R)$ and $K := {\rm SO_2}(\R)$, then we can identify ${\rm SL}_2(\R)\backslash (\bH \times \bH)$ with the double coset space $K\backslash G/K$, and hence $\eta$ can be seen as a Radon measure on this space. Secondly, one can show that there is a well-defined homeomorphism
\begin{equation}\label{iotaIntro}
\iota: K\backslash G/K \to [1, \infty), \quad KgK \mapsto \frac{1}{2}\trace(g^\top g),
\end{equation}
and hence $\eta$ corresponds to a Radon measure on $[1, \infty)$. We now offer several descriptions of this measure.

Firstly, we can identify $C_c(K\backslash G/K)$ with the convolution algebra $C_c(G, K)$ of bi-$K$-invariant functions on $G$ via pullback along the canonical projection $G \to K\backslash G/K$. From this identification $C_c(K\backslash G/K)$ inherits the structure of a $*$-algebra. Secondly, we can also identify every $f \in C_c(K\backslash G/K)$ with a radial function $f_{\bH^2}$ on the hyperbolic plane. We then obtain the following description of $\eta$.
\begin{proposition}[General formula for the auto-correlation measure]\label{AC1} The auto-correlation measure $\eta$ is the unique Radon measure on $K\backslash G/K$ such that for all $f \in C_c(K\backslash G/K)$,
\[
\eta(f \ast f^*) = \int_{\Omega_{(\Lambda, w)}} \left| \int_{\bH^2} f_{\bH^2} \, d\mu\right|^2 \, d\nu(\mu).
\]
\end{proposition}
Using results from \cite{BHP1} we obtain the following alternative description. Here, we denote by $\mathcal F \subset G \times G$ a fundamental domain for the $\Gamma$-action on $G \times G$ and by $m_G$ a suitably normalized choice of left-Haar measure on $G$. Moreover, given $f \in C_c(G, K)$ we denote by ${}_Kf_K \in C_c(K\backslash G/K)$ the function given by ${}_Kf_K(KgK) = f(g)$. \begin{theorem}[Auto-correlation formula for weighted model sets]\label{AC2} The auto-correlation measure $\eta$ is the unique Radon measure on $K\backslash G/K$ such that for all $f \in C_c(G, K)$,
\[
\eta({}_Kf_K \ast ({}_Kf_K)^*) \quad = \quad \int_{G} \int_G {\bf 1}_{\mathcal F}(g, h) \left|\sum_{(\gamma_1, \gamma_2) \in \Gamma} f(\gamma_1 g){\bf 1}_W(\gamma_2h)\right|^2\, dm_G(g)\, dm_G(h)
\]
where $\mathcal F$ is a fundamental domain for $\Gamma$ in $G \times G$. Equivalently, $\eta$ is the unique Radon measure on $K\backslash G/K$ such that for all $f \in C_c(G, K)$,
\[
\eta({}_Kf_K \ast ({}_Kf_K)^*)  \quad = \quad  \sum_{(\gamma_1, \gamma_2) \in \Gamma} (f \ast f^*)(\gamma_1) ({\bf 1}_W \ast {\bf 1}_{W^{-1}})(\gamma_2).
\]
\end{theorem}
Finally, denote by $(B_t)$ the ball of radius $t$ around $i$ in the hyperbolic plane, and define $F_t := \{g \in G\mid g.i \in B_t\} \subset G$.
\begin{theorem}[Sampling formula for the auto-correlation]\label{AC3}
The auto-correlation measure $\eta$ is the unique Radon measure on $K\backslash G/K$ such that for all $f \in C_c(G, K)$
\[
\eta({}_Kf_K) = \lim_{t\to \infty} \frac{1}{m_G(F_t)} \sum_{x \in \Lambda \cap F_t} \sum_{y \in \Lambda} f(xy^{-1}).
\]

\end{theorem}
From Proposition \ref{AC1} one sees that the auto-correlation measure is positive-definite on $K\backslash G/K$ in the sense that
\[
\eta(f \ast f^*)\geq 0 \quad \text{for all }f \in C_c(K\backslash G/K).
\]
However, if we consider $\eta$ as a Radon measure on $[1, \infty) \subset \R$ via the identification \eqref{iotaIntro}, then $\eta$ is not a positive definite Radon measure on $\R$ in these coordinates. We can remedy this by applying the so-called Harish transform and obtain an evenly positive-definite distribution on $\R$. To state the result, we denote by $C_c^\infty(\R)_{\rm ev} \subset C_c^\infty(\R)$ the subspace of
\emph{even functions}, i.e.\ functions satisfying $f(t) = f(-t)$. The dual space $\mathcal D(\R)_{\rm ev} :=C_c^\infty(\R)_{\rm ev}^*$ can be identified with the subspace of $\mathcal D(\R) = C_c^\infty(\R)^*$ consisting of those distributions which are invariant under the reflection at $0$, and hence we refer to elements of $\mathcal D(\R)_{\rm ev}$ as \emph{even distributions}. A distribution is \emph{positive-definite} if $\xi(\phi\ast \phi^*) \geq 0$ for all $\phi \in C_c^\infty(\R)$, and we call an even distribution $\xi$ \emph{evenly positive-definite} if $\xi(\phi \ast \phi^*) \geq 0$ for all $\phi \in C_c^\infty(\R)_{\rm ev}$. Given $\xi \in \mathcal D(\R)_{\rm ev}$ and $\phi \in C^\infty(\R)_{\rm ev}$ we write $\int_{0}^\infty \phi(t)\, d \xi(t) := \xi(\phi)$.
\begin{theorem}[Auto-correlation as a positive-definite distribution]\label{AC4}
If $\eta$ denotes the auto-correlation measure considered as a Radon measure on $[1, \infty)$, then the formula
\[
\xi(\phi) :=  \frac{-1}{2\pi} \int_1^\infty \int_{-\infty}^\infty \frac{\phi'(\arcosh(t+v^2/2))}{\sqrt{(t+v^2/2)^2-1}} \, dv\, d\eta(t) \quad (\phi \in C_c^\infty(\R)_{\rm ev})
\]
defines an evenly positive-definite distribution $\xi \in \mathcal D(\R)_{\rm ev}$, and for all $\psi \in C_c^\infty([1, \infty))$ we have
\[
\eta(\psi) =  \int_{0}^\infty \int_{-\infty}^\infty \psi(\cosh(t) + u^2/2)\,du\, d\xi(t).
\]
In particular, $\eta$ is uniquely determined by $\xi$.
\end{theorem}
In view of the theorem we refer to $\xi \in \mathcal D(\R)$ as the \emph{auto-correlation distribution} of $(\Lambda, w)$. In general, $\xi$ is not a tempered distribution, i.e.\ it does not extend to a continuous linear functional on the Schwartz space $\mathcal S(\R)$. The reason for this is that, unlike for model sets in the Euclidean plane, the number of elements of a model sets contained in a ball of radius $R$ in the hyperbolic plane grows exponentially (rather than polynomially) in $R$. While tempered distributions can be studied by their \emph{real} Fourier transform, to study non-tempered distributions one needs to employ a certain \emph{complex} Fourier transform. We will see in the sequel article \cite{BHP3}, that $\xi$ is uniquely determined by its complex Fourier transform, which we will show to be a pure-point positive Radon measure supported on a certain $1$-dimensional subset of $\C$. This should be seen as the analog of pure point diffraction of Euclidean quasi-crystals in the hyperbolic setting.

\subsection{The general picture}
All of the results mentioned this introduction work in much larger generality: The basic theory of auto-correlation measures up to Proposition \ref{AC1} can be developed for translation bounded measures with uniformly locally bounded orbit (in particular, weighted Delone sets of finite local complexity) in an arbitrary locally compact homogeneous metric space. The formula in Theorem \ref{AC2} still works for arbitrary weighted model sets in this very general setting. The approximation formula in Theorem \ref{AC3} is more restrictive. It works for example if $G$ is amenable and $(F_t)$ is a weakly admissible F\o lner sequence as define in \cite{BHP1}. It also works in many non-amenable situations, for example if $G$ satisfies the conclusion of the Howe--Moore theorem, $(G, K)$ is a Gelfand pair and $(F_t)$ is an arbitrary bi-$K$-invariant weakly admissible sequence. Notably this covers the case of balls in Riemannian symmetric space, in particular balls in the hyperbolic plane. Theorem \ref{AC4} makes use of an identification between the $*$-algebra $C_c^\infty(K\backslash G/K)$ and the subalgebra of $C_c^\infty(\R)$ consisting of even functions. This can be extended to semisimple Lie groups: If $G$ is a semisimple Lie group of real rank $n$ with maximal compact subgroup $K$, then by work of Harish-Chandra the $*$-algebra $C_c^\infty(K\backslash G/K)$ is isomorphic to the subalgebra of $C_c^\infty(\R^n)$ consisting of functions which are invariant under a certain finite reflection group, the so-called Weyl group of $G$. In this case, the auto-correlation measure can be identified with a Weyl group invariant positive-definite distribution on $\R^n$, and we will investigate complex Fourier transforms of these distributions in the sequel article \cite{BHP3}.

\subsection{Organization of the article} 
In Section \ref{DeloneGeneral} we discuss Delone sets in and translation bounded measures on a general proper lcsc metric space $X$. We explain how Delone sets in $X$ give rise to translation bounded measures (Proposition \ref{TBDirac}) and how they can be lifted to Delone sets in the isometry group provided $X$ is homogeneous (Corollary \ref{LiftDelone2}). As special cases we obtain Lemma \ref{IntroLift} and Proposition \ref{IntroLift2} from the introduction. We conclude by discussing weighted model sets as important examples.

In Section \ref{SecHull} we define the hull dynamical system of a translation bounded measure. Corollary \ref{CorHullCompact} yields a compactness result generalizing Proposition \ref{HullCompactIntro} from the introduction (as well as \cite[Thm.\ 4]{BaakeLenz} in the abelian case). Proposition \ref{HullVsHull} relates hulls of Delone sets (as studied in \cite{BHP1}) to the hull dynamical systems of their Dirac combs and Lemma \ref{HullNatural} yields naturality of these systems under proper equivariant maps. Together with results from \cite{BHP1} these imply (the general form of) Theorem \ref{UniqueErgodicityIntro} (cf.\ Corollary \ref{HullWeightedModelSet}).

Section \ref{SecAutocor} introduces the auto-correlation measure of a sufficiently nice translation bounded measure (including Dirac combs of weighted model sets). The most general form of Proposition \ref{AC1} is given in Proposition \ref{PropAutocorFormula}. We then specialize to the case of weighted model sets and derive the auto-correlation formulas from Theorem \ref{AC2} in Corollary \ref{Autocor1} and Proposition \ref{PropAutocor2}.

The final two sections are logically independent. Section \ref{SecSl2} discusses the notion of auto-correlation distribution (for ${\rm SL}_2(\R)$ and more general semisimple Lie groups) and establishes Theorem \ref{AC4}. Section \ref{SecApprox} establishes a sampling formula like the one in Theorem \ref{AC3} in large generality (Theorem \ref{ApproximationThm}).

The appendix collects various useful facts concerning convolution structures on double coset spaces, in particular concerning the existence of certain approximate identities.

\subsection{Notational conventions}
Throughout this article $G$ will always denote a unimodular lcsc group $G$ and the letter $K$ is reserved to denote a compact subgroup of $G$. We denote by $m_G$ a fixed choice of Haar measure on $G$ and by $m_K$ the Haar probability measure on $K$. Moreover, we will use the following notations.

\begin{remark}[Notations concerning point sets in groups] Given subsets $A, B \subset G$ we denote by $AB := \{ab\mid a \in A, b \in B\}$ the \emph{product set} of $A$ and $B$. Similarly, we define $A^{-1} := \{a^{-1} \mid a \in A\}$ and write $A^{n+1} := AA^n$ for iterated product sets. To avoid confusion with Cartesian products of sets we will usually write $X^{\times 2} := X \times X$ and $X^{\times (n+1)} := X \times X^{\times n}$ for iterated Cartesian products of a set $X$ with itself, except in standard notations like $\R^n$ or $\C^n$.
\end{remark}

\begin{remark}[Notations concerning $G$-spaces]
By a $G$-space we shall mean a lcsc space $X$ together with an action of $G$ on $X$ which is jointly continuous in the sense that the map $G \times X \to X$, $(g, x) \mapsto g.x$ is continuous. If there exists a metric $d$ on $X$ which defines the topology and is invariant under $G$ in the sense that $d(g.x, g.y) = d(x,y)$ for all $g \in G$ and $x, y \in \Omega$, then $(X, d)$ is called an \emph{isometric $G$-space}. If $\Omega$ is a compact $G$-space (in particular metrizable), then we sometimes call $\Omega$ a \emph{topological dynamical system} (TDS) over $G$.
\end{remark}

\begin{remark}[Notations concerning function spaces]
If $X$ is a lcsc space, then we denote by $C_c(X)$, $C_0(X)$ and $C_b(X)$ the function spaces of complex-valued compactly supported continuous functions, continuous functions vanishing at infinity and continuous bounded functions respectively.

If $(X, \nu)$ is a measure space and $f,g \in L^2(X, \nu)$, then we denote by \[\langle f, g \rangle_{X} := \langle f, g \rangle_{(X, \nu)} := \int_X f \cdot \overline{g} \, d\nu\]  the $L^2$-inner product. Contrary to the convention in \cite{BHP1} we will choose all our inner products to be anti-linear in the second variable.

Given a function $f: G \to \C$ we denote by $\bar f$, $\check f$ and $f^*$ respectively the functions on $G$ given by \[\bar f(g) := \overline{f(g)}, \quad \check f(g) := f(g^{-1}) \quad \text{and} \quad f^*(g) := \overline{f(g^{-1})}.\]
Given $f \in C_c(G)$ and $x,y \in G$ we define $L_xf(y) := f(x^{-1}y)$ and $R_xf(y) := f(yx)$. 
\end{remark}
\begin{remark}[Notations concerning measures] If $X$ is a lcsc space, then we denote by $M(X) = C_c(X)^*$ the Banach space of complex Radon measure on $X$. We write $M_b(X)$ for the subspace of finite complex measures (i.e. $\mu \in M(X)$ with $|\mu|(X)< \infty$) and $M^+(X)$ for the subset of (positive) Radon measures. Finally, we denote by $M_b^+(X)$ the space of bounded Radon measures on $X$ and by ${\rm Prob}(X) \subset M^+_b(X)$ the space of probability measures on $X$. We identify $\mu \in M(X)$ with the corresponding linear functional on $C_c(X)$ and write $\mu(f) := \int_X f \, d\mu$ for $f \in C_c(X)$.
\end{remark}

{\bf Acknowledgement.} The authors acknowledge financial support from the Young Excellence Grant (1142331) at Gothenburg Center of Advanced Studies (GoCas) and L\"angmanska kulturfonden BA19-1702. The authors are grateful for the hospitality shown to them by the Departments of Mathematics at Technion, Paderborn University, Giessen University and Chalmers University during the time that they worked on this paper. They thank Stephan Elsenhans for help in computing the explicit example in Subsection \ref{IntroWMS} and Stefan Witzel for  allowing us to use his beautiful picture of a Vornoi tiling associated with a weighted model set in the hyperbolic plane.

\section{Point sets and measures in proper homogeneous spaces}\label{DeloneGeneral}

\subsection{Metrics on proper homogeneous spaces}
Recall that $G$ denotes a unimodular lcsc group with Haar measure $m_G$ and that $K<G$ denotes a compact subgroup with Haar probability measure $m_K$. We denote by $K\backslash G$, $G/K$ and $K\backslash G/K$ respectively the quotients of $G$ by the left-action of $K$, right-action of $K$ and action of $K \times K$ respectively and denote by
\[
{}_Kp: G \to K\backslash G, \quad p_K: G \to G/K \qand {}_Kp_K: G \to K\backslash G/K
\]
the canonical projections. We will always topologize $K\backslash G$, $G/K$ and $K\backslash G/K$ with the quotient topology with respect to these projections, so that ${}_Kp$ and $p_K$ are open, closed and proper. While $K\backslash G$ and $G/K$ are homeomorphic, we will prefer to work with the left-quotient $K\backslash G$. In the sequel we will refer to $K\backslash G$ as a \emph{proper homogeneous space} of $G$.

\begin{example} Recall from the introduction that a metric space $(X,d_X)$ is called \emph{homogeneous} if its isometry group $G := {\rm Is}(X,d)$ acts transitively on $X$. If $X$ is a lcsc space (hence a proper metric space), then $G$ is a lcsc group with respect to the topology of pointwise convergence and for every $x_0 \in X$ the stabilizer $K:= {\rm Stab}_G(x_o)$ is compact \cite[Lemma 5.B.4]{CdlH}. Moreover, by the open mapping theorem, the map $K\backslash G \to X$, $Kg \mapsto g^{-1}(x_o)$ is a homeomorphism, hence every proper homogeneous metric space is a proper homogeneous space in our sense.
\end{example}

Our standing assumptions on $G$ imply that there exists a proper, continuous and right-invariant metric on $G$ which automatically defines the given topology on $G$ (Struble's theorem, see \cite[2.B.4]{CdlH}). We call any such metric \emph{right-admissible}. By averaging over $K$ we can produce a right-admissible metric on $G$ which is moreover left-$K$-invariant; we call such a metric \emph{$(K,G)$-admissible}. 

Given a $(K,G)$-admissible metric $d_G$, we can define a metric on $K\backslash G$ by setting
\[
d(Kg, Kh) := \min_{k \in K}d_G(kg, h).
\]
We refer to this metric as the \emph{induced metric} on $K\backslash G$. It is proper, continuous, $G$-invariant and induces the quotient topology on $K\backslash G$.

The group $G$ acts on $K\backslash G$ by $g.Kx := Kxg^{-1}$, and if $d$ is a metric on $K\backslash G$ which is induced by a $(K,G)$-admissible metric, then the action of $G$ on $(K\backslash G)$ is by isometries.

\begin{example} Let $(X,d_X)$ be a proper homogeneous metric space with isometry group $G$ with the topology of pointwise convergence. By \cite[Prop.\ 4.4.6]{Papadopoulos} a right-admissible metric on $G$ is given by
\[
d_G(g,h) := \sup_{x \in X}\, d_X(g^{-1}(x), h^{-1}(x)) e^{-d_X(x, x_0)}
\]
for any basepoint $x_0 \in X$. If we set $K := {\rm Stab}_G(x_0)$, then this metric is even $(K,G)$-admissible \cite[Prop.\ 4.4.4]{Papadopoulos}, and hence induces a metric $d$ on $K\backslash G$. Under the canonical identification $K\backslash G \cong X$ this metric is given by
\[
d(g^{-1}(x_0), h^{-1}(x_0)) = \min_{k \in K} \sup_{ x\in X} d_X(g^{-1} \circ k^{-1}(x), h^{-1}(x)) e^{-d_X(x,x_0)} \quad (g,h \in G).
\]
By definition we have for all $g, h \in G$,
\[
d(g^{-1}(x_0), h^{-1}(x_0)) \geq \min_{k \in K} d_X(g^{-1} \circ k^{-1}(x_0), h^{-1}(x_0)) e^{-d_X(x_0,x_0)} = d_X(g^{-1}(x_0), h^{-1}(x_0)).
\]
Note that if $f \in G$, then since $f$ is an isometry we have
\begin{eqnarray*}
d_X(f(x), x) e^{-d(x,x_0)} &\leq& (d_X(f(x_0), x_0) + d(x, x_0) + d(f(x), f(x_0)) e^{-d(x,x_0)}\\
&\leq& d_X(f(x_0), x_0) e^{-d(x,x_0)} + 2d(x,x_0)e^{-d(x, x_0)}.
\end{eqnarray*}
Applying this to $f := h \circ g^{-1} \circ k^{-1}$ and using $k^{-1}(x_0) = x_0$ and $\max \{2te^{-t} \mid t \geq 0\} = 2e^{-1}$ we obtain
\[
d(g^{-1}(x_0), h^{-1}(x_0)) \quad = \quad  \min_{k \in K} \sup_{ x\in X}  d_X(f(x), x) e^{-d(x,x_0)} \quad = \quad  d_X(g^{-1}(x_0), h^{-1}(x_0)) +2e^{-1}.
\]
\end{example}
We may thus record:
\begin{proposition}[Lifting metrics up to quasi-isometry]\label{ddXCompare} If $(X,d_X)$ is a proper homogeneous metric space with isometry group $G$ and point stablizer $K$, then there exists a metric $d$ on $X$ induced by a $(K,G)$-admissible metric $d_G$ on $G$ such that
\begin{equation}\label{ddXCompareineq}
\pushQED{\qed} 
d_X(x,y)\leq d(x,y) \leq d_X(x,y) + 2e^{-1} \quad (x,y \in X).
\qedhere
\popQED
\end{equation}
\end{proposition}
\subsection{Weighted Delone sets and translation bounded measures}
\begin{definition}[Terminology concerning weighted point sets]\label{RemPointSets}
Let $(X,d)$ be a lcsc metric space and let $\Lambda \subset X$ be a subset.
\begin{enumerate}
\item $\Lambda$ is called \emph{discrete} if every subset of $\Lambda$ is open in $\Lambda$ with respect to the subspace topology, and \emph{locally finite} if it is closed and discrete, or equivalently the intersection with every pre-compact subset of $X$ is finite. 
\item $\Lambda$ is called \emph{$r$-uniformly discrete} for some $r>0$ if for all $x,y \in X$ with $x \neq y$ we have $d(x,y) > r$. It is called \emph{$R$-relatively dense} for some $R>0$ if its $R$-neighbourhood $N_R(\Lambda)$ in $X$ coincides with $X$. It is called a \emph{$(r,R)$-Delone set} if it is $r$-uniformly discrete and $R$-relatively dense for some $R>r>0$. These notions depend on the choice of metric $d$. Any uniformly discrete set is locally finite. 
\item We denote by ${\rm LF}(X)$ the collection of all locally finite subsets of $X$, by ${\rm U}_r(X)$ the collection of all $r$-uniformly-discrete subsets of $X$, and by ${\rm Del}_{r,R}(X)$ the collection of all $(r,R)$-Delone subsets.
\item A subset $\mathcal A \subset {\rm LF}(X)$ is called \emph{uniformly locally finite} if for every pre-compact subset $K \subset X$ there exists a constant $C(K)$ such that $|\Lambda \cap K| \leq C(K)$ for all $\Lambda \in \mathcal A$.
\item If $\Lambda \subset X$ is a locally finite set, then we refer to a bounded function $w: \Lambda \to \C$ as a \emph{weight function} and to $(\Lambda, w)$ as a \emph{weighted point set}. It is called uniformly discrete or Delone if the underlying set $\Lambda$ is.
\end{enumerate}
\end{definition}
\begin{remark}[Weighted point sets as measures] If $(X,d)$ is a lcsc metric space, then every $\Lambda \in {\rm LF}(X)$ defines a Radon measure $\delta_\Lambda$ on $X$, called the \emph{associated Dirac comb}, by the formula
\[
\delta_\Lambda := \sum_{x \in \Lambda} \delta_x, \quad \text{i.e. } \delta_\Lambda(f) = \sum_{x \in \Lambda} f(x) \; \text{ for all } f\in C_c(X).
\]
Since $\Lambda \cap \supp(f)$ is finite for every $f \in C_c(X)$, this is well-defined. For $\Lambda = \emptyset$ we have $\delta_\emptyset = 0$ by convention. Similarly, if $(\Lambda, w)$ is a weighted point set, then we may define
\[
\delta_{(\Lambda,w)} := \sum_{x \in \Lambda} w(x) \cdot \delta_x \in M(X).
\]
In the sequel we will thus often think of locally finite sets as Radon measures via the embedding
\begin{equation}\label{LocallyFinitevsRadonMeasure}
{\rm LF}(X) \hookrightarrow M^+(X), \quad \Lambda \mapsto \delta_\Lambda,
\end{equation}
and similarly for weighted point sets.
\end{remark}
Dirac combs of (weighted) uniformly discrete point sets in isometric $G$-spaces have additional properties such as the following (cf. \cite{BaakeLenz}):
\begin{definition}[Translation bounded measures] Let $X$ be a lcsc $G$-space. A complex Radon measure $\mu \in M(X)$ is called \emph{translation bounded} with respect to $G$ (or \emph{$G$-bounded} for short) if for every compact $L \subset X$ we have
\[
\sup_{g \in G}|\mu(gL)| < \infty.
\]
We denote the space of $G$-bounded measures on $X$ by $\mathcal T_b(G \curvearrowright X)$.
\end{definition}
\begin{remark}[Translation bounded measures and uniformly locally finite orbits] If $X$ is a $G$-space and $\Lambda \subset X$ is a locally finite subset, then $\delta_\Lambda$ is $G$-bounded if and only if for every compact subset $K \subset X$ we have $\sup_{g \in G}|\Lambda \cap gK| < \infty$. This means precisely that the orbit $G.\Lambda$ of $\Lambda$ is uniformly locally finite.
\end{remark}
\begin{proposition}[Uniformly discrete subsets and translation bounded measures]\label{TBDirac} Let $(\Lambda, w)$ be a uniformly discrete weighted point set in an isometric lcsc $G$-space $(X,d)$. Then $\delta_{(\Lambda,w)}$ is $G$-bounded. In particular, if $\Lambda$ is uniformly discrete, then its orbit is uniformly locally finite.
\end{proposition}
\begin{proof} We fix a basepoint $o \in X$. Given an $r$-uniformly discrete set $\Lambda$ and $R>0$ we choose $x_1, \dots, x_n \in X$ such that $B_R(o)$ is covered by the balls $B_{r/2}(x_i)$. Then for every $g \in G$ we have
\[
g.B_R(o) \subset \bigcup_{i=1}^n gB_{r/2}(x_i) = \bigcup_{i=1}^n B_{r/2}(gx_i),
\]
and since every ball of radius $r/2$ contains at most one element of $\Lambda$ we have $\delta_{(\Lambda,w)}(g.B_R(o)) \leq n \cdot \|w\|_\infty$, whence $\delta_{(\Lambda,w)}$ is translation bounded.
\end{proof}
The converse implication is not true, even for $X = G$ and the constant weight $1$. A counterexample is given by the subset $\Lambda := \Z \cup \{n + \frac1 n \mid n \in \Z\{0\}\}$ of $\R$, which is not uniformly discrete, but whose $\R$-orbit is uniformly locally finite and whose Dirac comb is therefore $\R$-bounded.

\subsection{Weighted Delone sets in groups and proper homogeneous spaces}
\begin{remark}[Delone sets in groups] As explained in \cite{BH}, all right-admissible metrics on $G$ define the same notion of uniformly discrete, respectively relatively dense subsets, and we refer to such subsets simply as uniformly discrete, respectively relatively dense subsets of $G$. Explicitly, $\Lambda \subset G$ is uniformly discrete if $e$ is not an accumulation point of $\Lambda\Lambda^{-1}$ and relatively dense if $L \Lambda = G$ for some compact subset $L \subset G$. In particular, $\Lambda \subset G$ is uniformly discrete if $\Lambda\Lambda^{-1}$ is locally finite, in which case $\Lambda$ is said to have \emph{(right-)finite local complexity} (FLC).
\end{remark}
The following analogous result holds in the case of proper homogeneous spaces:
\begin{lemma}[Metric independence of Delone sets]\label{DeloneKG}
Let $\Lambda \subset K\backslash G$, denote by ${}_Kp: G \to K\backslash G$ the canonical projection, set $\Xi := {}_Kp^{-1}(\Lambda)$ and let $d_G$ be $(K,G)$-admissible with induced metric $d$.
\begin{enumerate}[(i)]
\item $\Lambda$ is relatively dense with respect to $d$ iff $\Xi$ is relatively dense in $G$.
\item $\Lambda$ is uniformly discrete with respect to $d$ iff there exists an identity neighbourhood $U$ in $G$ such that
$\Xi\Xi^{-1}\cap UK \subset K$.
\end{enumerate}
In particular, these notions are independent of the choice of $d$.
\end{lemma}
\begin{proof} (i) Assume that for every $Kg \in K\backslash G$ there exists $\lambda = K\xi \in \Lambda$ with $d(Kg, \lambda) < R$. Then $d_G(g, \xi) < R + {\rm diam}(K)$, hence $\Xi$ is relatively dense in $G$. The converse is immediate.
(ii) $\Lambda$ is $r$-uniformly discrete if and only if for all $\xi, \xi' \in \Xi$ we have either $K\xi = K\xi'$ or $d(K\xi, K\xi')>r$. Equivalently, for all $g = \xi\xi'^{-1} \in \Xi\Xi^{-1}$ we have either $g \in K$ or $d_G(k,g) = d_G(k\xi, \xi')>r$ for all $k \in K$. By right-invariance of the metric this means that $B_r(e)(\Xi\Xi^{-1} \setminus K) \cap K = \emptyset$, i.e.\ $\Xi\Xi^{-1} \cap B_r(e)K \subset K$.
\end{proof}
In view of the lemma we will refer to a uniformly discrete, relative dense or Delone set in $K\backslash G$ with respect to some (hence any) metric $d$ induced from a $(K,G)$-admissible metric simply as a  uniformly discrete, relative dense or Delone set in $K\backslash G$. The lemma motivates the following definition:
\begin{definition} A subset $\Xi \subset G$ is called \emph{$K$-uniformly discrete} if $\Xi\Xi^{-1}\cap UK \subset K
$ for some identity neighbourhood $U$ in $G$.
\end{definition}
In this terminology a subset $\Xi \subset G$ is uniformly discrete iff it is $\{e\}$-uniformly discrete. Note that in general a $K$-uniformly discrete set need not be uniformly discrete.
\begin{proposition}[Lifting Delone sets] \label{DeloneSet}
For a subset $\Lambda \subset K\backslash G$ the following are equivalent:
\begin{enumerate}[(i)]
\item $\Lambda$ is a Delone set.
\item There exists a $K$-uniformly discrete relatively dense set $\Xi \subset G$ such that ${}_Kp(\Xi) = \Lambda$
\item There exists a $K$-uniformly discrete Delone set $\Delta \subset G$ such that ${}_Kp(\Delta) = \Lambda$.
\end{enumerate}
\end{proposition}
\begin{proof} (i) $\Rightarrow$ (ii) By Lemma \ref{DeloneKG} we may choose $\Xi := {}_Kp^{-1}(\Lambda)$.

(ii) $\Rightarrow$ (iii) Let $\Delta \subset \Xi$ be a subset, which intersects each ${}_Kp$-fiber of $\Xi$ in a single point and fix a right-admissible metric $d_G$ on $G$. If $\Xi$ is $R$-relatively dense for some $R>0$ with respect to $d_G$, then $\Delta$ is $(R+{\rm diam}(K))$-relatively dense, and $\Delta$ is $K$-uniformly discrete as a subset of $\Xi$. Finally, if $\delta_1, \delta_2 \in \Delta$ are distinct, then $\delta_1\delta_2^{-1} \not \in K$, and hence $d(\delta, \delta_2) = d(\delta_1\delta_2^{-1}, e) \geq {\rm dist}(\Xi\Xi^{-1} \setminus K, K) > 0$, which shows that $\Delta$ is uniformly discrete. 

(iii) $\Rightarrow$ (i) Since $\Xi := {}_Kp^{-1}(\Lambda)$ contains $\Delta$, it is relatively dense, and hence $\Lambda$ is relatively dense by Lemma \ref{DeloneKG}. If $\lambda_1$ and $\lambda_2$ are two distinct points in $\Lambda$ with respective pre-images $\delta_1$ and $\delta_2$ in $\Delta$, then 
\[
d(\lambda_1, \lambda_2) = \min_{k \in K} d_G(k\delta_1, \delta_2) =  \min_{k \in K} d_G(\delta_2\delta_1^{-1}, k) \geq {\rm dist}(\Xi\Xi^{-1} \setminus K, K) > 0,
\]
hence $\Lambda$ is uniformly discrete.
\end{proof}
\begin{remark} From Proposition \ref{DeloneSet} one deduces the following:
\begin{enumerate}
\item Every projection of a $K$-uniformly discrete Delone set is again Delone.
\item Every Delone set in $K\backslash G$ is a projection of a Delone set in $G$. By analogy, we will say that
a subset of $K\backslash G$ has FLC if it is the projection of an FLC set in $G$.
\item The proof of Lemma \ref{DeloneKG} shows that $\Xi \subset G$ is $K$-uniformly discrete if and only if for every right-admissible metric ${\rm dist}(\Xi\Xi^{-1} \setminus K, K) > 0$. This implies that every FLC set in $G$ is $K$-uniformly discrete, and hence every FLC set in $K\backslash G$ is uniformly discrete.
\end{enumerate}
We conclude in particular, that if $\Lambda \subset G$ is an FLC Delone set, then ${}_Kp(\Lambda)$ is a Delone set
\end{remark}
\begin{remark}[Application to proper homogeneous metric spaces] Let $(X, d_X)$ be a proper homogeneous lcsc metric space with isometry group $G$ and point stabilizer $K$ and let $d$ be a metric as in Proposition \ref{ddXCompare}. We claim that a subset $\Lambda \subset X$ is Delone with respect to $d_X$ if and only if it is Delone with respect to $d$. Indeed, that relative denseness carries over is immediate from \eqref{ddXCompareineq}, and this inequality also implies that any $d_X$-uniformly discrete subset is $d$-uniformly discrete with the same constant. Now assume that $\Lambda$ is $r$-uniformly discrete with respect to $d$, but not with respect to $d_X$. Then there exist elements $x_n \neq y_n$ in $\Lambda$ with $d_X(x_n, y_n) \leq \frac{1}{n}$. Choose $g_n \in G$ such that $g_nx_n = x_0$; then $d_X(x_0, g_n(y_n)) \leq \frac{1}{n}$ and thus $g_n(y_n) \to x_0$. Since $d$ is continuous this implies 
\[
d(x_n, y_n) = d(g_n(x_n), g_n(y_n)) = d(x_0, g_n(y_n)) \to d(x_0, x_0) = 0,
\]
hence for sufficiently large $n$ we have $d(x_n, y_n) < r$, which is a contradiction.
\end{remark}
Together with Proposition \ref{DeloneSet} we deduce:
\begin{corollary}[Lifting Delone sets, II]\label{LiftDelone2} Let $(X, d_X)$ be a lcsc proper metric space. Assume that $G<{\rm Is}(X)$ acts transitively on $X$ with point-stabilizer $K$. Then a subset $\Lambda \subset X$ is Delone if and only if it is the orbit of a $K$-uniformly discrete Delone set in $G$.\qed
\end{corollary}

\begin{remark}[Push-forwards of point sets and measures]\label{RemPushForward} 
Let ${}_Kp: G \to K\backslash G$ denote the canonical projection.
\begin{enumerate}
\item We have seen that if $\Lambda \subset G$ is uniformly discrete, then the \emph{naive push-forward} ${}_Kp(\Lambda)$ need not be uniformly discrete.
\item On the other hand, since  ${}_Kp$ is proper, it induces a push-forward map ${}_Kp_*: M(G) \to M(K\backslash G)$, and the push-forward of a translation bounded measure will always be translation bounded, i.e. ${}_Kp_*$ restricts to a map ${}_Kp_*: M(G) \to M(K\backslash G)$. In this sense translation bounded measures behave better under push-forward than uniformly discrete sets.
\item  Let $\Lambda \subset G$ be uniformly discrete. Then the \emph{weighted push forward} ${}_Kp_*\Lambda$ of $\Lambda$ is the weighted point set $({}_Kp(\Lambda), w)$ with weight function given by $w(x) = |\Lambda \cap {}_Kp^{-1}(x)|$. The definition is made in such a way that
\[
\delta_{{}_Kp_*\Lambda} = {}_Kp_*\delta_\Lambda.
\]
This implies in paticular that $\delta_{{}_Kp_*\Lambda}$ is translation bounded, and hence $\Lambda$ has a a uniformly locally finite $G$-orbit.
\item If $\widetilde{w}$ is a weight function on $\Lambda$, then we can similarly define ${}_Kp_*(\Lambda, \widetilde{w})$ as the unique weighted point set such that
\[
\delta_{{}_Kp_*(\Lambda, \widetilde{w})} = {}_Kp_*\delta_{(\Lambda,\widetilde{w})}.
\]
Explicitly, the weight function $w$ of ${}_Kp_*(\Lambda, \widetilde{w})$ is given by
\[
{w}(x) = \sum_{y \in {}_Kp^{-1}(x)} \widetilde{w}(y).
\] 
\end{enumerate}
Note that if $\Lambda \subset G$ is an FLC Delone set, then ${}_Kp(\Lambda)$ is a Delone set, and hence ${}_Kp_*\Lambda$ is a weighted Delone set.
\end{remark}
\subsection{Weighted model sets in proper homogeneous spaces}
We now discuss the class of examples which motivated the current series of papers. In \cite{BHP1} we introduced the notion of a model set in a lcsc group $G$. To state the definition, we recall that a discrete subgroup $\Gamma$ of a lcsc group $G$ is a \emph{lattice} if $\Gamma\backslash G$ admits a $G$-invariant probability measure, and a \emph{uniform lattice} if $\Gamma\backslash G$ is moreover compact. 
\begin{definition}
\begin{enumerate}
\item  A \emph{cut-and-project-scheme} is a triple $(G, H, \Gamma)$ where $G$ and $H$ are lcsc groups and $\Gamma < G \times H$ is a lattice which projects injectively to $G$ and densely to $H$. A cut-and-project scheme is called \emph{uniform} if $\Gamma$ is moreover a uniform lattice.
\item If $(G, H, \Gamma)$ is a cut-and-project scheme and $p_G : G\times H \to G$ denotes the projection onto the first coordinate, then for every pre-compact\footnote{In \cite{BHP1} we insisted that $W$ be compact, but this is actually never used and precludes some naturally arising examples. (DOUBLE CHECK THIS!)} set $W \subset G$ the subset 
\[
\Lambda(G,H, \Gamma, W) := p_G(\Gamma \cap (G \times W)) \subset G
\]
is called a \emph{weak model set} in $G$ with \emph{window} $W$. A weak model set is called uniform if $\Gamma$ is a uniform lattice.
\item A weak (uniform) model set is called a \emph{(uniform) model set} if its window $W$ has non-empty interior. It is called a \emph{regular model set} if $W$ is Jordan-measurable with dense interior, ${\rm Stab}_H(W) = \{e\}$ and $\partial \overline{W} \cap p_H(\Gamma) = \emptyset$, where $p_H :G \times H \to H$ is the projection onto the second factor.
\end{enumerate}
\end{definition}
Now let $\Lambda \subset G$ be a model set, let $K < G$ be a compact subgroup and ${}_Kp: G \to K\backslash G$ denote the canonical projection. Then the weighted push-forward ${}_Kp_*\Lambda$ is a weighted Delone set of the form $({}_Kp(\Lambda), w_\Lambda)$.
\begin{definition}\label{DefWeightedModelSet} The weighted Delone set ${}_Kp_*\Lambda = ({}_Kp(\Lambda), w_\Lambda)$ is called a \emph{weighted model set} in $K\backslash G$, and $w_\Lambda$ is called its \emph{canonical weight function}. $(\pi(\Lambda), w_\Lambda)$ is called a \emph{regular} (respectively \emph{uniform}) weighted model set, if $\Lambda$ has the corresponding property.
\end{definition}
We emphasize that in our terminology, a weighted model set is not just a weighted Delone set with underlying Delone set ${}_Kp(\Lambda)$. Rather, when we speak of weighted model sets, we always assume that the weight function is the canonical one.
\begin{remark} Let $\Lambda$ be a model set in $G$ and ${}_Kp_*\Lambda = ({}_Kp(\Lambda), w)$. If the weight function $w$ is trivial (i.e. $w(x) = 1$ for all $x \in {}_Kp(\Lambda)$), then we may identify ${}_Kp_*\Lambda$ with the underlying set ${}_Kp(\Lambda)$. We then call ${}_Kp(\Lambda)$ simply a \emph{model set} in $K\backslash G$. By definition, ${}_Kp(\Lambda)$ is a model set if and only if $\Lambda\Lambda^{-1} \cap K = \{e\}$.
\end{remark}

\begin{convention}\label{regular} Throughout this article, all weighted model sets (uniform or not) are assumed to be regular.
\end{convention}

\section{Hulls of Delone sets and translation bounded measures}\label{SecHull}

\subsection{Topologies on point sets and measures}
\begin{remark}\label{BasicTopologies} Let $(X,d)$ be a lcsc metric space. We define the following topologies on the various sets of point sets and measures defined above.
\begin{enumerate}
\item We equip $M(X)$ with the weak-$*$-topology with respect to $C_c(X)$. This topology is second-countable, since $X$ is, and the action of $G$ on $M(X)$ is jointly continuous with respect to the weak-$*$-topology. Indeed, if  $(g_n, \mu_n) \to (g, \mu)$ in $G \times M(X)$, then for every $\varphi \in C_c(X)$ there exists a compact set $K$ containing the supports of all of the functions $g_n.\varphi$. If we set $C := {\sup} \mu_n(K)$, then $C < \infty$ and
\begin{eqnarray*}
|{\mu_n}(g_n.\varphi) - \mu(g.\varphi)| &\leq& |{\mu_n}(g_n.\varphi) - {\mu_n}(g.\varphi)|+|{\mu_n}(g.\varphi) - \mu(g.\varphi)|\\
&\leq& C \cdot \|g_n.\varphi -g.\varphi\|_\infty+|{\mu_n}(g.\varphi) - \mu(g.\varphi)|.
\end{eqnarray*}
Since both terms converge to $0$ as $n \to \infty$, this shows $g_n.\mu_n \to g.\mu$.

\item The set $\mathcal C(X)$ of all closed subsets of $X$ carries a natural compact metrizable topology called the \emph{Chabauty--Fell topology}. A sequence $(P_n)$ in $\mathcal C(X)$ converges to $P$ with respect to this topology if and only if the following two properties hold:
\begin{enumerate}[(CF1)]
\item If $(n_k)$ is an unbounded sequence of natural numbers and $p_{n_k} \in P_{n_k}$ converge to $p \in X$, then $p \in P$.
\item For every $p \in P$ there exist elements $p_n \in P_n$ such that $p_n \to p$.
\end{enumerate}
From this characterization one sees immediately that the $G$ action on $\mathcal C(X)$ is jointly continuous with respect to the Chabauty--Fell topology.
\item For every $R>r>0$ we have inclusions ${\rm Del}_{r,R}(X) \subset {\rm U}_r(X) \subset {\rm LF}(X) \subset \mathcal C(X)$ and we denote by $\tau_{CF}$ the restrictions of the Chabauty-Fell topology to either of these spaces.
\item We obtain another topology on ${\rm LF}(X)$ and its subspaces by pulling back the weak-$*$-topology on $M(X)$ via the embedding ${\rm LF}(X) \hookrightarrow M(X)$ from \eqref{LocallyFinitevsRadonMeasure}. We refer to this topology as the \emph{measure topology} on ${\rm LF}(X)$. The measure topology on weighted point sets is defined accordingly.
\end{enumerate}
\end{remark}
\begin{proposition}[Chabauty--Fell topology vs.\ measure topology]\label{CFvsMeasureTopology} Let $(X,d)$ be a lcsc metric space. Then for every $r > 0$ the Chabauty--Fell topology and the measure topology coincide on the subset ${\rm U}_r(X) \subset {\rm LF}(X)$ and define a compact metrizable topology on ${\rm U}_r(X)$.
\end{proposition}
Explicitly this means that a sequence of $r$-uniformly discrete sets $\Lambda_n$ converges to $\Lambda$ in the Chabauty--Fell topology if and only if for every $f \in C_c(X)$ we have
\[
\sum_{x_n \in \Lambda_n} f(x_n) \to \sum_{x \in \Lambda} f(x). 
\] 
For the proof we follow roughly the same strategy as in the abelian case, cf.\ \cite[Thm.\ 4]{BaakeLenz}. We need the following lemma:
\begin{lemma}[Convergence of uniformly discrete sets]\label{UrConv}
Suppose that $\Lambda_n \to \Lambda$ in ${\rm U}_r(X)$ with respect to the Chabauty--Fell topology. Then for every $p \in \Lambda$ and $0<\epsilon<r/2$ there exists $n_0(p, \epsilon) \in \mathbb N$ such that for all $n \geq n_0(p, \epsilon)$, 
\[
|B_{\epsilon}(p) \cap \Lambda_n| = 1.
\]
\end{lemma}
\begin{proof} If $p \in \Lambda$, then by (CF2) there exist $p_n \in \Lambda_n$ with $p_n \to p$, hence for every $\epsilon >0$ there exists $n_0(p)$ such that for all $n \geq n_0$ we have $p_n \in B_{\epsilon}(p) \cap \Lambda_n$. If $\epsilon < r/2$, then $p_n$ is necessarily unique by the triangle inequality.
\end{proof}
\begin{proof}[Proof of Proposition \ref{CFvsMeasureTopology}] The space $({\rm U}_r(X), \tau_{CF})$ is compact and second countable (see e.g.\ \cite{BH}). It thus suffices to show that the embedding $({\rm U}_r(X), \tau_{CF}) \hookrightarrow M(X)$, $\Lambda \mapsto \delta_\Lambda$ is (sequentially) continuous.

Thus assume that $\Lambda_n \to \Lambda$ in $({\rm U}_r(X), \tau_{CF})$ and let $f \in C_c(X)$. Let $K_0 := {\rm supp}(f)$, and let $K_1$ be a compact set containing a $10r$-neighbourhood of $K_0$.
%$K_1 \subset K_2\subset X$ compact such that $K_1$ contains a $10r$-neighbourhood of $K_0$ and $K_2$ contains a $100r$-neighbourhood of $K_0$. 
Note that $K_1 \cap \Lambda$ is finite, say $K_1 \cap \Lambda = \{p^{(1)}, \dots, p^{(N)}\}$. Now fix $\epsilon \in (0, r/2]$ and with the notation of Lemma \ref{UrConv} set
\[
n_0 := \max\{n_0(p^{(j)}, \epsilon)\mid j \in \{1, \dots, N\}\}.
\]
Then for all $n \geq n_0$ we have $B_\epsilon(p^{{j}}) \cap \Lambda_n = \{p_n^{(j)}\}$. We claim that for all but finitely many $n \geq n_0$ we have
\[
\Lambda_n \cap K_0 \subset \{p_n^{(j)} \mid j \in \{1, \dots, N\}\}.
\]
Indeed, otherwise we would find an unbounded sequence $n_k$ and elements $p_{n_k} \in \Lambda_{n_k} \cap K_0$ such that $d(p_{n_k}, \Lambda \cap K_1) > \epsilon$. Passing to a further subsequence we may assume that $p_{n_k}$ converges to some $p \in X$, and then $p \in \Lambda$ by (CF1). Since $p$ is not contained in $K_1$, it has distance at least $10r$ from $K_0$. But then at most finitely many of the $p_{n_k}$ can be contained in $K_0$, which is a contradiction. We deduce that for sufficiently large $n$ we have
\begin{eqnarray*}
\left| {\delta_{\Lambda_n}}(f) -{\delta_\Lambda}(f)\right|&=& \left|\sum_{j=1}^{N} f(p_n^{(j)}) -\sum_{j=1}^N f(p^{(j)})\right|\\
&\leq& \sum_{j=1}^{N} |f(p_n^{(j)}) - f(p^{(j)})|\\
&\leq& N \cdot \sup\{|f(x)-f(y)|\mid x,y \in N_r(K_1), d(x,y) < \epsilon\}. 
\end{eqnarray*}
Since $f$ is uniformly continuous on compacta, we deduce that ${\delta_{\Lambda_n}}(f) \to {\delta_\Lambda}(f)$, which finishes the proof.
\end{proof}
\subsection{The hull of a translation bounded measure}
From now on let $G$ be a lcsc group and let $(X,d)$ be an isometric $G$-space. We fix a basepoint $o \in X$. 
\begin{remark}[$C$-translation bounded measures]
Since every compact subset of $X$ is contained in some ball around $o$, a complex measure $\mu \in M(X)$ is translation bounded if and only if for every $R>0$ there exists $C(R) > 0$ such that for all $g \in G$,
\[
\mu(gB_R(o)) \leq C(R).
\]
If this holds for a fixed function $C: (0, \infty) \to (0, \infty)$, then we say that $\mu$ is \emph{$C$-translation bounded} (with respect to $d$ and $o$). Given a function $C: (0, \infty) \to (0, \infty)$ we then denote by $\mathcal T_{C}(G \curvearrowright X) \subset \mathcal T_b(G \curvearrowright X)$ the collection of all complex Radon measures on $X$ which are $C$-translation bounded with respect to $d$ and $o$. By definition, the set $\mathcal T_{C}(G \curvearrowright X)$ is invariant under the natural action of $G$ on complex Radon measures.
\end{remark}
\begin{lemma}[Compactness properties of translation bounded measures] For every function $C: (0, \infty) \to (0, \infty)$ the subspace $\mathcal T_{C}(G \curvearrowright X)  \subset M(X)$ is a compact metrizable space, and the $G$-action on $\mathcal T_{C}(G \curvearrowright X)$ is continuous.
\end{lemma}
\begin{proof} By Remark \ref{BasicTopologies}, $M(X)$ is second countable (hence metrizable) and the action of $G$ on $M(X)$ is jointly continuous. It thus remains to show that $\mathcal T_{C}(G \curvearrowright X)\subset M(X)$ is compact. Since $M(X)$ carries the subspace topology with respect to the embedding 
\begin{equation}\label{Weak*Top}
M(X) \; \hookrightarrow \prod_{\varphi \in C_c(X)} \R, \quad \mu \mapsto (\mu(\varphi)),
\end{equation}
where the right-hand side is given the product topology, this amounts to showing that $\mathcal T_{C}(G \curvearrowright X)$ has compact image under this embedding. If $\varphi \in C_c(X)$ and $\mu \in \mathcal T_{C}(G \curvearrowright X)$, then 
\[
|\mu(\varphi)| \leq \|\varphi\|_\infty \cdot \mu({\rm supp}(\varphi)),
\]
and hence the image of $\mathcal T_{C}(G \curvearrowright X)$ under the embedding is bounded in each coordinate. Moreover, for every $\mu \in M(X)$ we have
\[
\mu(gB_R(o)) = \sup\{|\mu(\varphi)| \mid \varphi \in C_c(X), {\rm supp}(\varphi) \subset gB_R(o), \|\varphi\|_\infty \leq 1\}.
\]
hence the conditions $\mu(gB_R(o)) \leq C(R)$ are closed conditions. This shows that $\mathcal T_{C}(G \curvearrowright X)$ is compact and finishes the proof.
\end{proof}
\begin{corollary}\label{CorHullCompact} If $\mu \in \mathcal T_b(G \curvearrowright X)$, then the orbit closure $\Omega_\mu := \overline{G.\mu}$ is a compact metrizable space and the action of $G$ on $\Omega_\mu$ is continuous.\qed
\end{corollary}
\begin{definition}[Hull of a translation-bounded measure] Given $\mu \in \mathcal T_b(G \curvearrowright X)$, the orbit closure
\[\Omega_\mu := \overline{G.\mu} \subset M(G)
\]
is called the \emph{hull} of $\mu$, and the TDS $G \curvearrowright \Omega_\mu$ is called the associated \emph{hull system}.
\end{definition}
\begin{remark}[Hulls systems of locally finite sets] In \cite{BHP1} we defined the notion of a hull of a locally finite subset of $G$. More generally, let $X$ be a lcsc $G$-space and let $\Lambda \in {\rm LF}(X)$. Then one defines the \emph{hull} of $\Lambda$ as the orbit closure $\Omega_\Lambda := \overline{G.\Lambda}$, where the orbit closure is taken with respect to the Chabauty--Fell topology. The associated TDS $G \curvearrowright \Omega_\Lambda$ is then called the \emph{hull system} of $\Lambda$.
\end{remark}
\begin{proposition}[Hulls of uniformly discrete sets vs.\ hulls of their Dirac combs]\label{HullVsHull} Let $(X,d)$ be an lcsc isometric $G$-space and let $\Lambda \subset X$ be uniformly discrete. Then there is a well-defined $G$-equivariant homeomorphism 
\[
\Omega_\Lambda \to \Omega_{\delta_\Lambda}, \quad \Lambda' \mapsto \delta_{\Lambda'}.
\]
\end{proposition}
\begin{proof} Assume that $\Lambda$ is $r$-uniformly discrete. Since $X$ is isometric, every translate $g\Lambda$ of $\Lambda$ is also $r$-uniformly discrete, and since $\mathcal U_r(X)$ is closed with respect to the Chabauty--Fell topology, we have $\Omega_\Lambda \subset \mathcal U_r(X)$. By Proposition \ref{CFvsMeasureTopology} the map $\Omega_\Lambda \hookrightarrow \mathcal U_r(X) \to M(X)$ given by $\Lambda' \mapsto \delta_{\Lambda'}$ is continuous, and it is evidently $G$-equivariant. It thus maps the orbit closure of $\Lambda$ onto the orbit closure of its image. By compactness, the resulting bijection is a homeomorphism.
\end{proof}
In this sense, the hull system of a uniformly discrete set can be seen as a special case of the hull system of a translation bounded measure. If $(\Lambda, w)$ is a weighted point set we define $\Omega_{(\Lambda, w)} := \Omega_{\delta_{\Lambda, w}}$ and refer to $\Omega_{(\Lambda, w)}$ as the hull system of $(\Lambda, w)$.

\begin{remark}[Punctured hull]
Given a translation bounded measure $\mu \in \mathcal T_b(G \curvearrowright X)$ we observe that it may happen that $0 \in \Omega_\mu$ even if $\mu \neq 0$. For example, if $\mu = \delta_\Lambda$ is the Dirac comb of a uniformly discrete set $\Lambda \subset G$, then this happens if and only if $\Lambda$ is not relatively dense (see \cite{BHP1}).  It is then convenient to remove the $G$-fixpoint and to consider also the \emph{punctured hull} $\Omega_\mu^\times := \Omega_\mu \setminus \{0\}$. By construction, $\Omega_\mu^\times$ is a lcsc space, and if $0 \in \Omega_\mu$, then $\Omega_\mu$ is its one-point-compactification. 
\end{remark}
\begin{lemma}[Naturality of the hull under push-forward]\label{HullNatural} Let $X$, $Y$ be lcsc $G$-spaces, let $\mu \in \mathcal T_b(G \curvearrowright X)$ and let $p: X \to Y$ be a proper continuous $G$-equivariant map. Then $p$ induces continuous $G$-factor maps
\[
p_*: \Omega_\mu \to \Omega_{p_*\mu}, \quad \mu' \mapsto p_*\mu' \qand p_*: \Omega^\times_\mu \to \Omega^\times_{p_*\mu}, \quad \mu' \mapsto p_*\mu'.
\]
\end{lemma}
\begin{proof} As in Remark \ref{RemPushForward} one observes that $p$ induces a continuous $G$-equivariant push-forward map $p_*: \mathcal T_b(G \curvearrowright X) \to \mathcal T_b(G \curvearrowright Y)$, since it is proper. It thus maps the orbit closure of $\mu$ surjectively onto the orbit closure of $p_*\mu$. This proves the first statement, and the second statement then follows from the fact that $0$ is the only pre-image of $0$ under $p_*$.
\end{proof}
As an application we can extend properties of hulls of model sets in lcsc groups from \cite{BHP1} to hulls of weighted model sets in proper homogeneous spaces. We remind the reader that by a weighted model set we always mean a projection of a model set together with its \emph{canonical} weight function (cf. Definition \ref{DefWeightedModelSet}). We also recall from Convention \ref{regular} that  all weighted model sets are implicitly assumed to be regular.

\begin{corollary}[Hulls of weighted model sets]\label{HullWeightedModelSet} The punctured hull of a weighted model set in a proper homogeneous $G$-space $K\backslash G$ is uniquely ergodic with respect to the $G$-action. If the weighted model set is uniform, then its hull is moreover minimal.
\end{corollary}
\begin{proof} Let $\Lambda \subset G$ be a model set and consider ${}_Kp_*\Lambda = ({}_Kp(\Lambda), w)$ in $K\backslash G$. By Lemma \ref{HullNatural} we have a continuous $G$-factor map ${}_Kp_*: \Omega^\times_\Lambda \to \Omega^\times_{{}_Kp_*\Lambda}$. 
Let $\mu$ be a probability measure on $G$ which is absolutely continuous with respect to Haar measure and whose support generates $G$ as a semigroup. It was established in \cite[Thm. 3.4]{BHP1} that there exists a unique $\mu$-stationary probability measure $\nu$ on $\Omega_\Lambda$, which is moreover $G$-invariant. In particular, ${}_Kp_*\nu$ defines a $G$-invariant probability measure on $\Omega^\times_{{}_Kp_*\Lambda}$. Now if $\nu'$ is any $G$-invariant probability measure on $\Omega^\times_{{}_Kp_*\Lambda}$, then its fiber
\[
\{\widetilde{\nu'} \in {\rm Prob}(\Omega_\Lambda)\mid {}_Kp_*\widetilde{\nu'} = \nu'\}
\]
is a compact convex $G$-invariant set, hence contains a fixpoint under the convolution action by $\mu$. By uniqueness, this fixpoint must coincide with $\nu$, and hence $\nu' = {}_Kp_*\nu$, showing that $\Omega^\times_{{}_Kp_*\Lambda}$ is uniquely ergodic. If $\Lambda$ is uniform, then $\Omega^\times_\Lambda$ is minimal by \cite[Prop. 3.3]{BHP1}, and this property descends to the continuous factor $\Omega^\times_{{}_Kp_*\Lambda}$.
\end{proof}
\begin{remark}[Invariant measure on the hull of a weighted model set]\label{InvariantMeasureWeightedModelSet} We record for later reference that the unique $G$-invariant probability measure on $ \Omega^\times_{{}_Kp_*\Lambda}$ is the push-forward of the unique $G$-invariant probability measure on $\Omega_\Lambda$ under the map ${}_Kp_*$.
\end{remark}

\subsection{Uniform local boundedness of the hull}
In the next section we are going to define a periodization map for the hull of a translation-bounded measure. Continuity of this map will depend on a property of the hull called uniform local boundedness, which we investigate in the present subsection.

Recall from Remark \ref{RemPointSets} that, given a lcsc space $X$, a subset $\Omega \subset {\rm LF}(X)$ is called d \emph{uniformly locally finite}, if for every pre-compact subset $K \subset X$ there exists a constant $C_K>0$ such that $|\Lambda \cap K| \leq C_K$ for all $\Lambda \in \Omega$. Similarly, a subset $\mathcal A \subset M(X)$ will be called \emph{uniformly locally bounded} if for every pre-compact subset $K \subset G$ there exists a constant $C_K>0$ such that $\mu(K) \leq C_K$ for all $\mu \in \mathcal A$. Here we are interested in conditions on translation bounded measures which guarantee uniform local boundedness of the hull.

We first consider the situation for point sets $\Lambda$ in an isometric lcsc $G$-space $(X,d)$. By Proposition \ref{TBDirac}, uniform discreteness of $\Lambda$ is enough to ensure uniform local finiteness of the orbit $G.\Lambda$, but it is not enough to ensure uniform local finiteness of the orbit \emph{closure} $\Omega_\Lambda$. In the case where $X = G$, the latter is implied by finite local complexity:
\begin{proposition}[Uniform local finiteness of FLC hulls]\label{FLCHull} Let $G$ be a lcsc group and $\Lambda \subset G$ be an FLC set. Then the hull $\Omega_\Lambda$ is uniformly locally finite, and more generally the hull $\Omega_{({\Lambda, w})}$ is uniformly locally bounded for every weight function $w : \Lambda \to [0, \infty)$.
\end{proposition}
The proof of the proposition makes use of the following lemma from \cite{BH}. Since we worked with left-FLC rather than right-FLC sets there, we repeat the short proof:
\begin{lemma}\label{P-1P} Let $G$ be a lcsc group and let $\Lambda \subset G$ be a locally finite subset. Then for all $\Lambda' \in \Omega_\Lambda$ we have $\Lambda'(\Lambda')^{-1} \subset \overline{\Lambda\Lambda^{-1}}$. 
\end{lemma}
\begin{proof} If $\Lambda' \in \Omega_\Lambda$ then there exist $g_n \in G$ such that $\Lambda g_n^{-1} \to \Lambda'$. By (CF2) we thus find for every $p,q \in \Lambda'$ sequences $(p_n)$ $(q_n)$ in $\Lambda$ such that $p_ng_n^{-1} \to p$ and $q_ng_n^{-1} \to q$. By continuity of multiplication and inversion in $G$ we obtain $p_nq_n^{-1} \to pq^{-1}$ and thus $\Lambda' (\Lambda')^{-1}\subset \overline{\Lambda\Lambda^{-1}}$.
\end{proof}
We will apply this in the following form:
. \begin{corollary}\label{FLCHull2} Let $\Lambda \subset G$ be an FLC set. Then there exists an open identity neighbourhood $U\subset G$ such that $|\Lambda' \cap Ug^{-1}| \leq 1$ for all $\Lambda' \in \Omega_\Lambda$ and $g \in G$.
\end{corollary}
\begin{proof} Since $\Lambda\Lambda^{-1}$  is locally finite, there exists an open identity neighbourhood $V$ such that $\Lambda\Lambda^{-1} \cap V = \{e\}$. Since $\Lambda\Lambda^{-1}$ is closed, Lemma \ref{P-1P} then shows that $\Lambda'(\Lambda')^{-1} \cap V = \{e\}$ for all $\Lambda' \in \Omega_\Lambda$. Now let $U \subset G$ be an open identity neighbourhood with $UU^{-1} \subset V$. Given $g \in G$ and $\Lambda' \in \Omega_\Lambda$ we either have $\Lambda' \cap Ug^{-1} = \emptyset$ or there exist $p \in \Lambda'$ and $u \in U$ such that $p = ug^{-1}$, i.e.\ $g = p^{-1}u$. In the latter case we have
\[
\Lambda' \cap Ug^{-1} = \Lambda'\cap Uu^{-1}p = (\Lambda'p^{-1} \cap Uu^{-1})p \subset (\Lambda'(\Lambda')^{-1} \cap V)p = \{p\},
\] 
hence $|\Lambda' \cap Ug^{-1}| \leq 1$ in either case.
\end{proof}
\begin{proof}[Proof of Proposition \ref{FLCHull}] Let $K \subset G$ be pre-compact. Choose $U$ as in Lemma \ref{FLCHull2} and cover $K$ by finitely many translates of $U$, say $K \subset Ug_1^{-1} \cup \dots \cup Ug_N^{-1}$. Then by the lemma we have for every $(\Lambda',w') \in \Omega_{\Lambda, w}$,
\[
\delta_{\Lambda', w'}(K) \leq \|w'\|_\infty \sum_{j=1}^N |\Lambda' \cap Ug_j^{-1}| \leq N \cdot \|w\|_\infty,
\]
and this bound is independent of $(\Lambda', w')$.
\end{proof}
In particular, Proposition \ref{FLCHull} implies that the hull of every model set is uniformly locally finite. To extend this result to weighted model sets we are going to use:
\begin{proposition}\label{ULBPushForward} Let $X$, $Y$ be lcsc $G$-spaces, let $\mu \in \mathcal T_b(G \curvearrowright X)$ and let $p: X \to Y$ be a proper continuous $G$-equivariant map. If $\mu$ has a uniformly locally bounded hull, then also $p_*\mu$ has a uniformly locally bounded hull.
\end{proposition}
\begin{proof} Let $L \subset Y$ be compact. Since $M := p^{-1}(L)$ is compact, there exists $C_M>0$ such that for all $\mu' \in \Omega_\mu$ we have $\mu(M) < C_M$. It follows that if $\mu' \in \Omega_{\mu}$, then
\[
p_*\mu'(L) = \mu'(p^{-1}(L)) = \mu(M) < C_M.
\]
Since by Lemma \ref{HullNatural} the map $p_*: \Omega_\mu \to \Omega_{p_*\mu}$ is onto, the proposition follows.
\end{proof}
Applying this to the canonical projection ${}_Kp: G \to K\backslash G$ and using Proposition \ref{FLCHull} we have arrived at the following result:
\begin{corollary}[Hulls of weighted model sets are uniformly locally bounded]\label{FLCHullDownstairs} If $\Lambda \subset G$ is an FLC set and $w$ is an arbitrary weight function on $\Lambda$, then the hull $\Omega_{{}_Kp_*(\Lambda, w)}$ is uniformly locally bounded. In particular, the hull of every weighted model set is uniformly locally bounded.\qed
\end{corollary}

\section{Auto-correlation measures for translation-bounded measures}\label{SecAutocor}

\subsection{The periodization map}
Let $X$ be a lcsc $G$-space. The following construction generalizes the periodization map of a locally finite set in $G$ as introduced in \cite{BH} to the case of a translation bounded measure in $X$.
\begin{proposition}[Periodization over the hull]\label{PropPeriodization} Let $\mu \in \mathcal T_b(G \curvearrowright X)$.
\begin{enumerate}[(i)]
\item For every $f \in C_c(X)$ the function
\[
\mathcal P_\mu f: \Omega_\mu \to \C, \quad \mathcal P_\mu f(\mu') := \int_X f d\mu'
\]
is well-defined and continuous on $\Omega_\mu$.
\item If $0 \in \Omega_\mu$, then $\mathcal P_\mu f(0) = 0$ for all $f\in C_c(G)$, hence we obtain a map $\mathcal P_\mu: C_c(X) \to C_0(\Omega^\times_\mu)$.
\item The map $\mathcal P_\mu$ is $G$-equivariant.
\end{enumerate}
\end{proposition}\begin{proof} The map $\mathcal P_\mu f$ is just the restriction of the evaluation map $M(X) \to \C$, $\mu' \mapsto \mu'(f)$, which is well-defined and continuous by definition of the weak-$*$-topology. This shows (i), and (ii) and (iii) are immediate from the definitions.
\end{proof}
\begin{definition} For $\mu \in \mathcal T_b(G \curvearrowright X)$ the $G$-equivariant map 
\[
\mathcal P_\mu: C_c(X) \to C_0(\Omega^\times_\mu), \quad \mathcal P_\mu f(\mu') := \int_X f \, d\mu'
\]
is called the \emph{periodization} map associated with $\mu$. 
\end{definition}
Continuity of the periodization map depends on uniform local boundedness of the hull: We equip $C_c(X)$ with the direct limit topology
\[
C_c(X) = \lim_{\to} (C(L), \|\cdot\|_\infty),
\]
where $L$ runs through all compact subsets of $X$. Then continuity of the periodization map $\mathcal P_\mu$ with respect to this topology means that for every compact subset $L \subset X$ there exists $C_L > 0$ such that if $f \in C_c(G)$ with ${\rm supp}(f) \subset L$, then
\begin{equation}
\|\mathcal P_\mu f\|_\infty \quad \leq \quad C_L \cdot \|f\|_\infty.
\end{equation}
\begin{proposition}[Continuity of the periodization map]\label{PContinuous} Let $\mu \in \mathcal T_b(G \curvearrowright X)$. If $\Omega_\mu$ is uniformly locally bounded, then the periodization map $\mathcal P_\mu: C_c(X) \to C(\Omega_\mu)$ is continuous with respect to the Frech\'et topology on $C_c(X)$.
\end{proposition}
\begin{proof} Let $L\subset X$ be compact. Since $\Omega_\mu$ is uniformly locally bounded there exists $C_{L} >0$ such that for all $\mu' \in \Omega_\mu$ we have $\mu'(L) < C_L$. But then for all $f \in C_c(G)$ with ${\rm supp}(f) \subset L$ we have
\[
\|\mathcal P_\mu f\|_\infty = \sup_{\mu' \in \Omega_\mu}\left| \mathcal P_\mu f(\mu')\right| =  \sup_{\mu' \in \Omega_\mu}\left| \int_X f \, d\mu'\right|  \leq  \mu'({\rm supp}(f)) \cdot \|f\|_\infty \leq C_L \cdot \|f\|_\infty.
\]
This shows that $\mathcal P_\mu$ is continuous and finishes the proof.
\end{proof}
\begin{remark}[Conditions ensuring continuity of the periodization map] The periodization map is continuous in each of the following cases:
\begin{enumerate}[(i)]
\item $X = G$ (with the action given by $g.x := xg^{-1}$) and $\mu = \delta_\Gamma$, where $\Gamma < G$ is a discrete subgroup. For this case it was established in \cite{BHS} that $\Omega_\Gamma^\times = \Gamma\backslash G$, and the periodization map $\mathcal P_\Gamma: C_c(G) \to C_0(\Gamma\backslash G)$ is given by the classical formula
\[
\mathcal P_\Gamma f(\Gamma g) = \sum_{\gamma \in \Gamma} f(\gamma g).
\]
\item $X=G$ (with the action given by $g.x := xg^{-1}$) and $\mu = \delta_\Lambda$, where $\Lambda$ is a (right-)FLC set in $G$. In this case, under the canonical identification of $\Omega^\times_\mu$ with $\Omega_\Lambda$, the periodization map is given by the formula
\[
\mathcal P_\Lambda: C_c(X) \to C_0(\Omega_\Lambda^\times), \quad \mathcal P_\Lambda f(\Lambda') := \sum_{x \in \Lambda'} f(x).
\]
Up to the issue of left- vs. right-action, this is precisely our definition of the periodization map from \cite{BH, BHP1}, hence the current definition is compatible with our previous one in this case.
\item $X=G$ and $\mu = \delta_{(\Lambda, w)}$, where $(\Lambda, w)$ is a weighted FLC set in $G$ (cf. Proposition \ref{FLCHull}).
\item $X= K\backslash G$ is a proper $G$-space and $\Lambda$ is an FLC set in $X$ and $w$ an arbitrary (bounded) weight function (cf. Corollary \ref{FLCHullDownstairs}).
\end{enumerate}
In particular, we have a continuous periodization map both for model sets in a lcsc group $G$ (as already defined in \cite{BHP1}) and weighted model sets in proper homgeoneous $G$-space. These are the two classes of examples in which we are most interested in the current article.
\end{remark}

If $\Gamma$ is a discrete subgroup of $G$, then by \cite[Sec.\ 1]{Raghunathan} the image of the periodization map $\mathcal P_\Gamma: C_c(G) \to C_0(\Gamma\backslash G)$ is given by $C_c(\Gamma\backslash G)$, hence a dense subalgebra of $C_0(\Gamma\backslash G)$. In the setting at hand we have the following weaker statement.
\begin{proposition}[Almost surjectivity of the periodization map]\label{PeriodizationAlmostSurjective}Let $\mu \in \mathcal T_b(G \curvearrowright X)$ be a translation bounded measure with uniformly locally bounded hull.
\begin{enumerate}[(i)]
\item If $\nu \in \Omega_\mu^\times$, then there exists $f \in C_c(X)$ such that $\mathcal P_\mu f(\nu) \neq 0$.
\item If $\nu_1, \nu_2 \in \Omega_\mu$ and $\nu_1 \neq \nu_2$, then there exists $f \in C_c(X)$ with $\mathcal P_\mu f(\nu_1) \neq \mathcal P_\mu f(\nu_2)$.
\item The algebra generated by the constant function $1$ and the image of the periodization map is dense in $C(\Omega_\mu)$.
\item The algebra generated by the image of the periodization map is dense in $C_0(\Omega_\mu^\times)$.
\end{enumerate}
\end{proposition}
\begin{proof} (i) and (ii) follow from the fact that a complex Radon measure is uniquely determined by the associated integral. %Consequently, if $\nu$ is a non-zero complex Radon measure on $X$ there exists $f \in C_c(X)$ with $\nu(f) \neq 0$ and if $\nu_1, \nu_2$ are distinct Radon measures on $X$ then there exists  $f \in C_c(X)$ with ${\nu_1}(f) \neq {\nu_2}(f)$. 
From (ii) we deduce that if $\nu_1, \nu_2 \in \Omega_\mu$ are distinct, then there exists $f \in C_c(X)$ such that $\mathcal P_\mu f(\nu_1) \neq \mathcal P_\mu f(\nu_2)$, and hence
\[
(\mathcal P_\mu f - {\nu_1}(f) \cdot 1)(\nu_1) = 0, \quad (\mathcal P_\mu f - {\nu_1}(f))(\nu_2) = \mathcal P_\mu f(\nu_2) - \mathcal P_\mu f(\nu_1) \neq 0.
\]
Thus the algebra generated by $1$ and the image of the periodization map separates points in $\Omega_\mu$, hence is dense in $C(\Omega_\mu)$ by the Stone--Weierstrass theorem. This shows (iii), and (iv) follows from (iii).
\end{proof}

\subsection{Functoriality of periodization}
\begin{remark}[Functoriality for measures]\label{PeriodizationFunctorial}
Let $X, Y$ be lcsc $G$-spaces, let $p: X \to Y$ be a proper continuous $G$-equivariant map and let $\mu \in \mathcal T_b(G \curvearrowright X)$. Then for all $f \in C_c(Y)$ and $g \in G$ we have
\[
(\mathcal P_{\mu} \circ p^*)f(g_*\mu) = \int_X p^*f \, dg_*\mu
=  \int_Y f \, dp_*(g_*\mu)
= \mathcal P_{p_*\mu} f(p_*(g_*\mu)) = (p^* \circ \mathcal P_{p_*\mu})f(g_*\mu),
\]
and hence the following diagram commutes:
\[\begin{xy}\xymatrix{
C_c(X) \ar[rr]^{\mathcal P_\mu} &&  C_0(\Omega_\mu^\times)\\
C_c(Y) \ar[rr]^{\mathcal P_{p_*\mu}} \ar[u]^{p^*}&& C_0(\Omega_{{p_*\mu}}^\times)\ar[u]_{(p_*)^*}
}\end{xy}\]
Most notably, this applies to the case, where $X = G$ and $Y = K\backslash G$ is a proper homogeneous $G$-space, and we will apply this in the case of (weighted) model sets below.
\end{remark}
\begin{caveat}[Failure of functoriality for point sets] If, in the situation of the previous remark, $\Lambda$ is a uniformly discrete subset of $X$ and $p(\Lambda)$ is the naive push-forward of $\Lambda$ to Y, then the following diagram need \emph{not} commute:
\[\begin{xy}\xymatrix{
C_c(X) \ar[rr]^{\mathcal P_\Lambda} &&  C_0(\Omega_\Lambda^\times)\\
C_c(Y) \ar[rr]^{\mathcal P_{p(\Lambda)}} \ar[u]^{p^*}&& C_0(\Omega_{{p(\Lambda)}}^\times)\ar[u]_{(p_*)^*}.
}\end{xy}\]
In this sense, the only ``natural'' push-forward of an FLC set $\Lambda \subset G$ to $K\backslash G$ via the canonical projection ${}_Kp: G \to K\backslash G$ is the weighted push-forward ${}_Kp_*\Lambda$, whereas the naive push-forward ${}_Kp(\Lambda)$ is not natural. This lack of functoriality for the naive push-forward is what led us to consider translation bounded measures on, rather than point sets in, $K\backslash G$.
\end{caveat}

\subsection{Correlation measures and the auto-correlation}
Throughout this subsection let $X$ be a lcsc $G$-space and $\mu\in \mathcal T_b(G \curvearrowright X)$ be a translation bounded measure with uniformly locally bounded punctured hull $\Omega^\times_\mu$.

Note that if $\nu \in M(\Omega^\times_\mu)$, then for every $n \in \mathbb N$ the linear functional on $C_c(X) \otimes \dots \otimes C_c(X)\subset C_c(X^{\times n})$ given by 
\[
f_1\otimes \dots \otimes f_n \mapsto \int_{\Omega_\mu^\times} \mathcal P_\mu f_1(\xi) \mathcal P_\mu f_2(\xi) \cdots \mathcal P_\mu f_n(\xi) \, d\nu(\xi)
\]
is continuous by Proposition \ref{PContinuous}. It is thus given by integration against a Radon measure on $X^{\times n}$, and we define:
\begin{definition} Given $n \in \mathbb N$ the measure $\eta_\nu^{(n)}  \in M(X^{\times n})$ given by
\begin{equation}\label{nthCorrelationMeasure}
\int_{X^{\times n}} f_1(x_1) \cdots f_n(x_n) d\eta_{\nu}^{(n)}(x_1, \dots, x_n) = \int_{\Omega_\mu^\times} \mathcal P_\mu f_1(\xi) \mathcal P_\mu f_2(\xi) \cdots \mathcal P_\mu f_n(\xi) \, d\nu(\xi).
\end{equation}
is called the \emph{$n$th correlation measure} of $\nu$.
\end{definition}
From Proposition \ref{PeriodizationAlmostSurjective}.(iv) we deduce:
\begin{corollary}[Correlation measures determine the measure] Every probability measure $\nu$ on $\Omega^\times_\mu$ is uniquely determined by the sequence $\eta_\nu^{(1)}, \eta_\nu^{(2)}, \dots$ of its correlation measures.\qed
\end{corollary}
It is an interesting question whether \emph{finitely} many correlation measures are enough to determine an arbitrary (i.e. not necessarily $G$-invariant) probability measure on $\Omega^\times_\mu$, but we will not pursue this here. Rather we focus on the case of a $G$-invariant Radon measure $\nu$ on $\Omega_\mu^\times$. The following observation is immediate from $G$-equivariance of the periodization map:
\begin{proposition}[Correlations of invariant measures are invariant]
If $\nu$ is a $G$-invariant Radon measure on $\Omega_\mu^\times$, then for every $n \in \mathbb N$ the correlation measure $\eta_\nu^{(n)}  \in M(X^{\times n})$ is $G$-invariant.\qed
\end{proposition}
As explained in Lemma \ref{WeilFormula} in the appendix, under suitable assumptions on $X$ and $G$ a $G$-invariant measure on $X^{\times n}$ corresponds to a measure on the quotient $G\backslash X^{\times n}$. Specifically for the second correlation measure we obtain:
\begin{corollary}[Auto-correlation measure]\label{CorAutoCor} Let $G$ be a unimodular lcsc group, $X$ be a $G$-space such that $X \times X$ is strongly proper as a $G$-space, and let $\mu\in \mathcal T_b(G \curvearrowright X)$ be a translation bounded measure with uniformly locally bounded hull. Then for every $G$-invariant Radon measure $\nu$ on $\Omega_\mu^\times$ there exists a unique Radon measure $\eta \in M(G \backslash(X \times X))$ such that for all $f \in C_c(X \times X)$,
\[
\pushQED{\qed} 
\int_{X \times X} f(x_1, x_2) d\eta_\nu^{(2)}(x_1, x_2) =\int_{G\backslash (X \times X)}  \left(\int_G f(gx_1, gx_2) dm_G(g)\right) d{\eta}(G(x_1, x_2)).
\qedhere
\popQED
\]
\end{corollary}
\begin{definition} In the situation of Corollary \ref{CorAutoCor} the measure $\eta\in M(G \backslash(X \times X))$ is called the \emph{auto-correlation measure} of the measure $\nu$.
\end{definition}
Higher auto-correlation measures can be defined similarly, but we will not consider them in the current article.

\subsection{Auto-correlation measures in the case of proper homogeneous $G$-spaces}
In this subsection we specialize to the case of a proper homogeneous $G$-space $X = K\backslash G$. We denote by ${}_Kp: G \to K\backslash G$ and ${}_Kp_K: G \to K\backslash G /K$ the canonical projections. Moreover, $\mu$ denotes a translation bounded measure on $X := K\backslash G$ with uniformly locally bounded punctured hull $\Omega^\times_\mu$ and $\nu$ denotes a $G$-invariant probability measure on $\Omega^\times_\mu$. We are interested in describing the auto-correlation measure $\eta$ of $\nu$ more explicitly.

We first consider the group case, where $K = \{e\}$ and $X=G$ with the $G$-action given by $g.x := xg^{-1}$. In this case, the quotient space $G\backslash (G \times G)$ can be identified with $G$ via the map $(g_1, g_2) \mapsto g_1g_2^{-1}$ and hence we obtain an isomorphism $M^+(G\backslash (G \times G)) \to M^+(G)$. In particular, we can consider the auto-correlation measure as a Radon measure on $G$ itself.
In the case where $\mu$ was a Dirac comb of an FLC set, this measure was described in \cite{BHP1}, and the result carries over to our setting at hand as follows:
\begin{proposition}[Autocorrelation formula, group case]\label{AutocorGroup} Let $\mu\in \mathcal T_b(G \curvearrowright G)$ be a translation bounded measure with uniformly locally bounded hull under the right-action of $G$ on itself and let $\nu$ be a $G$-invariant probability measure on $\Omega_\mu^\times$. Then the auto-correlation measure $\eta \in M^+(G)$ is the unique Radon measure such that
\begin{equation}\label{AutocorrelationGroupCase}
\eta(f \ast f^*) = \|\mathcal P_\mu f\|^2_{L^2(\Omega_\mu^\times, \nu)} \quad \text{for all }f \in C_c(G).
\end{equation}
\end{proposition}
\begin{proof} Let $f \in C_c(G)$. On the one hand we observe that
\[
\eta^{(2)}(f \otimes \overline{f}) = \int_{\Omega_\mu^\times} \mathcal P_\mu f \cdot \mathcal P_\mu \overline f \, d\nu =  \int_{\Omega_\mu^\times} \mathcal P_\mu f \cdot  \overline{\mathcal P_\mu f}\, d\nu= \|\mathcal P_\mu f\|^2_{L^2(\Omega_\mu^\times, \nu)}.
\]
On the other hand, in view of the identification $G\backslash (G \times G) \cong G$ above the measure $\eta \in M^+(G)$ satisfies
\[
\eta^{(2)}(f \otimes \overline{f}) = \eta(F),
\]
where for $g_1, g_2 \in G$ we have
\[
F(g_1g_2^{-1}) = \int_G (f \otimes \overline{f})(g_1g, g_2g) \, dm_G(g) = \int f(g)\overline{f}(g_2g_1^{-1}g) \, dm_G(g) = f\ast f^*(g_1g_2^{-1}).
\]
We conclude that $\eta(f \ast f^*) = \eta(F) = \eta^{(2)}(f \otimes \overline{f})$, which yields \eqref{AutocorrelationGroupCase}.

Finally, $\eta$ is uniquely determined by Formula \eqref{AutocorrelationGroupCase} in view of Corollary \ref{ff*Dense}.
\end{proof}
If $f,g \in C_c(G)$, then it follows from \eqref{AutocorrelationGroupCase} and polarization that
\[
\eta(f \ast g^*) = \langle \mathcal P_\mu f, \mathcal P_\mu g\rangle_{L^2(\Omega_\mu^\times, \nu)}.
\]
Thus if $(\rho_n)$ denotes a convenient approximate identity in $C_c(G)$ as in Remark \ref{ConvenientApproximateIdentity} consisting of real-valued symmetric functions with supports contained in a fixed compact set, then
\[
\eta(f) = \lim_{n \to \infty}  \langle \mathcal P_\mu f, \mathcal P_\mu\rho_n \rangle_{L^2(\Omega_\mu^\times, \nu)}.
\]
We now want to given a similar characterization of the autocorrelation measure for a general proper homogeneous space $K\backslash G$.
Observe first that if $X = K\backslash G$, then we have a natural identification
\[
G \backslash (X \times X) \to K\backslash G/K, \quad [(Kg_1, Kg_2)] \mapsto Kg_1g_2^{-1}K.
\]
We may thus consider the auto-correlation measure $\eta$ as a Radon measure on the double coset space $K\backslash G/K$, and hence as a linear functional on $C_c(K\backslash G/K)$. We recall from Subsection \ref{DoubleCosetConvolution} in the appendix that $C_c(K\backslash G/K)$ carries a natural convolution structure and involution such that ${}_Kp_K^*: C_c(K\backslash G/K) \to C_c(G, K)$ becomes an isomorphism of $*$-algebras. 

The canonical projection ${}_{(K)}p_K: K\backslash G \to K\backslash G/K$ induces an embedding ${}_{(K)}p_K^*: C_c(K\backslash G/K) \hookrightarrow C_c(K\backslash G)$ and by abuse of notation we denote the composition
\[
C_c(K\backslash G/K) \xrightarrow{{}_{(K)}p_K^*} C_c(K \backslash G) \xrightarrow{\mathcal P_\mu} C_0(\Omega^\times_\mu).
\]
also by $\mathcal P_\mu$. Explicitly, this means that
\[
\mathcal P_\mu f(\mu') := \mu'({}_{(K)}p_K^*f) \quad (f \in C_c(K\backslash G/K), \mu' \in \Omega_\mu^\times).
\]
With this abuse of notation understood we have:
\begin{proposition}[Autocorrelation formula, general case]\label{PropAutocorFormula} Let $\mu\in \mathcal T_b(G \curvearrowright K\backslash G)$ be a translation bounded measure with uniformly locally bounded hull and let $\nu$ be a $G$-invariant probability measure on $\Omega_\mu^\times$. Then the auto-correlation measure $\eta \in M^+(K\backslash G/K)$ is the unique Radon measure such that
\begin{equation}\label{AutocorrelationGeneralCase}
\eta(f \ast f^*) = \|\mathcal P_\mu f\|^2_{L^2(\Omega_\mu^\times, \nu)} \quad \text{for all }f \in C_c(K\backslash G/K).
\end{equation}
\end{proposition}
As in the group case, we can use a polarization argument to give a formula for the measure $\eta \in M^+(K\backslash G/K)$: If $\rho_n$ is a convenient approximate identity in $C_c(K\backslash G/K)$ (cf.\ Remark \ref{ConvenientApproximateIdentity}), then
\[
\eta(f) = \lim_{n \to \infty} \langle \mathcal P_\mu f, \mathcal P_\mu \rho_n \rangle_{L^2(\Omega_\mu^\times, \nu)} \quad \text{for all }f \in C_c(K\backslash G/K).
\]
\begin{proof}[Proof of Proposition \ref{PropAutocorFormula}] Let  $f \in C_c(K\backslash G/K)$ and denote by ${}_{(K)}p_K: K\backslash G \to K\backslash G/K$ the canonical projection. We set $f_\dagger := {}_{(K)}p_K^*f$ so that, by definition, $\mathcal P_\mu f = \mathcal P_\mu f_\dagger$. Then
\[
\eta^{(2)}(f_\dagger \otimes \overline{f_\dagger})=  \int_{\Omega_\mu^\times} \mathcal P_\mu f_\dagger \cdot  \overline{\mathcal P_\mu f_\dagger}\, d\nu= \|\mathcal P_\mu f_\dagger\|^2_{L^2(\Omega_\mu^\times, \nu)} = \|\mathcal P_\mu{f}\|^2_{L^2(\Omega_\mu^\times, \nu)} .
\]
On the other hand, in view of the identification $G\backslash (X \times X) \cong K\backslash G/K$ above the measure $\eta \in M^+(K\backslash G/K)$ satisfies
\[
\eta^{(2)}(f_\dagger \otimes \overline{f_\dagger}) = \eta(F),
\]
where for $g_1, g_2 \in G$ we have
\begin{eqnarray*}
F(Kg_1g_2^{-1}K) &=& \int_G (f_\dagger \otimes \overline{f_\dagger})(Kg_1g, Kg_2g) \, dm_G(g) \quad = \quad \int f(Kg_1gK)\overline{f}(Kg_2gK) \, dm_G(g)\\
&=& \int_G {}_Kp_K^*f(g_1g) {}_Kp_K^*f(g_2g) \, dm_G(g)
%&=& \int f(g)\overline{f}(g_2g_1^{-1}g) \, dm_G(g) = f\ast f^*(g_1g_2^{-1}).
\end{eqnarray*}
The same computation as in the proof of Proposition \ref{AutocorGroup} (but applied to ${}_Kp_K^*f$ instead of $f$) now yields
\[
F(Kg_1g_2^{-1}K) = {}_Kp_K^*f \ast ({}_Kp_K^*f)^*(g_1g_2^{-1}) = {}_Kp_K^*(f \ast f^*)(g_1g_2^{-1}) =  f \ast f^*(Kg_1g_2^{-1}K).
\]
This shows that $F = f \ast f^*$ and thus $\eta(f \ast f^*) = \eta(F) = \eta^{(2)}(f_\dagger \otimes \overline{f_\dagger})$, which yields \eqref{AutocorrelationGeneralCase}.

As in the group case, $\eta$ is uniquely determined by Formula \eqref{AutocorrelationGeneralCase} in view of Corollary \ref{ff*Dense}.
\end{proof}
In Subsection \ref{DoubleCosetConvolution} in the appendix we also define a function $f \ast f^* \in C_c(K\backslash G/K)$ for a given function $f \in C_c(K\backslash G)$. With this definition understood, the above proof actually shows that \eqref{AutocorrelationGeneralCase} holds for all $f \in C_c(K\backslash G)$. 
\begin{remark}[Positive-definiteness of the auto-correlation measure]\label{AutocorPosDef}
In the group case, the auto-correlation measure $\eta \in M^+(G)$ is ``positive-definite'' or of ``positive type'' in the usual sense, i.e.\ $\eta(f \ast f^*) \geq 0$ for all $f \in C_c(G)$. In the case of a proper homogeneous space, the auto-correlation measure $\eta \in M^+(K\backslash G/K)$ has the analogous positivity property that $\eta(f \ast f^*) \geq 0$ for all $f \in C_c(K\backslash G/K)$ for the given $*$-algebra structure on $C_c(K\backslash G/K)$. In this sense, the auto-correlation measure is positive-definite also in the case of proper homogeneous spaces.
\end{remark}

\subsection{Naturality of correlation and auto-correlation measures}
In this section we consider a translation bounded measure $\mu$ on $G$ (with respect to the right-action); we will assume that the punctured hull $\Omega^\times_\mu$ is uniformly locally bounded and that there exists a $G$-invariant probability measure $\nu$ on $\Omega^\times_\mu$. We also consider the push-forward ${}_Kp_*\mu$ of $\mu$ under the canonical projection ${}_Kp: G \to K\backslash G$.

By Proposition \ref{ULBPushForward}, the punctured hull $\Omega_{{}_Kp_*\mu}^\times$ is also uniformly locally bounded, and by Lemma \ref{HullNatural} the map ${}_Kp$ induces a continuous $G$-factor map ${}_Kp_*: \Omega_\mu \to \Omega_{{}_Kp_*\mu}$. It follows from $G$-equivariance that the push-forward ${}_Kp_*\nu$ of $\nu$ with respect to this map is a $G$-invariant probability measure on $\Omega_{{}_Kp_*\mu}$. 

We are going to compare the correlation measures $\widehat{\eta}^{(n)}$ of $\nu$ to the correlation measures $\eta^{(n)}$ of ${}_Kp_*\nu$ and the auto-correlation measure $\widehat{\eta}$ of $\nu$. to the correlation measure $\eta$ of ${}_Kp_*\nu$ and by $\eta$ the auto-correlation measure of ${}_Kp_*\nu$. For this we denote by ${}_Kp^{\times n}: G^{\times n} \to (K\backslash G)^{\times n}$ and ${}_Kp_K: G \to K\backslash G /K$ the canonical projections.

\begin{proposition}[Naturality of the auto-correlation]\label{AutocorNatural} The correlation measures and the auto-correlation measure of $\nu$ and ${}_Kp_*\nu$ are related by the formulas
\[
\eta^{(n)} = {}_Kp^{\times n}_*(\widehat{\eta}^{(n)}) \in M^+((K\backslash G)^n) \qand \eta =({}_Kp_K)_*(\widehat{\eta}) \in M^+(K\backslash G/K).
\]
\end{proposition}
\begin{proof} The second statement follows from the first one (for $n=2$) since the diagram
\[\begin{xy}\xymatrix{
M(G\times G)^G \ar[r] \ar[d]_{ {}_Kp^{\times 2}_*}& M(G\backslash (G \times G))\ar[r]&M(G)\ar[d]^{{}_Kp_K^*}\\
M(X\times X)^G \ar[r] & M(G\backslash (X\times X))\ar[r]&M(K\backslash G/K)
}\end{xy}\]
commutes. The first statement follows in turn from Remark \ref{PeriodizationFunctorial}, since for all $f_1, \dots, f_n \in C_c(X)$ we have
\begin{eqnarray*}
{{}_Kp^{\times n}_*}\eta^{(n)}(f_1 \otimes \dots \otimes f_n) &=& \int_{G^{\times n}} {}_Kp^*f_1(g_1) \cdots {}_Kp^*f_n(g_n)d\eta^{(n)}_{\nu}(g_1, \dots, g_n)\\
&=& \int_{\Omega_\mu^\times} \mathcal P_\mu {}_Kp^*f_1(\xi) \mathcal P_\mu {}_Kp^*f_2(\xi) \cdots \mathcal P_\mu {}_Kp^*f_n(\xi) \, d{\nu}(\xi)\\
&=& \int_{\Omega_\mu^\times} {}_Kp^*\mathcal P_{{}_Kp_*\mu}f_1(\xi) {}_Kp^*\mathcal P_{{}_Kp_*\mu}f_2(\xi) \cdots  {}_Kp^*\mathcal P_{{}_Kp_*\mu} f_n(\xi) \, d{\nu}(\xi)\\
&=& \int_{X^{\times n}}f_1(x_1) \cdots f_n(x_n) \, d\eta^{(n)}_{{}_Kp_*\nu}(x_1, \dots, x_n)\\
&=& \pushQED{\qedhere} 
\eta^{(n)}(f_1 \otimes \dots \otimes f_n).\qedhere \popQED
\end{eqnarray*}
The proposition follows.
\end{proof}

\subsection{Auto-correlation measures of weighted model sets}
We have seen in \cite{BHP1} that we can obtain an explicit formula for the auto-correlation measure of a model set $\Lambda$ in $G$. Using functoriality, we now extend this result to weighted model sets in proper homogeneous spaces.

Since the results in \cite{BHP1} were stated for the left-action of $G$ on itself, and since we prefer to work with the right-action here, we briefly restate the relevant results in the group case in our current notation. We fix a cut-and-project scheme $(G, H, \Gamma)$ and consider a model set of the form $\Lambda = \Lambda(G, H, \Gamma, W)$ with window $W \subset H$. We recall our convention that model sets are assumed to be regular in this article. We also recall that the punctured hull $\Omega_\Lambda^\times$ of $\Lambda$ is uniformly locally finite and uniquely ergodic and denote by $\widehat{\nu}$ the unique $G$-invariant probability measure on $\Omega_\Lambda^\times$.

By our standing assumptions on $W$, the characteristic function ${\bf 1}_W$ is a compactly-supported  bounded measurable function on $H$. Since $\Gamma$ is a uniformly discrete subset of $G\times H$ we can define the periodization of an arbitrary compactly-supported measurable functions $F: G \times H \to \C$ by the same formula as in the continuous case, i.e.\
\[
\mathcal P_\Gamma F(\Gamma(g,h)) := \sum_{(\gamma_1, \gamma_2) \in \Gamma} f(\gamma_1g, \gamma_2h).
\]
With this notation, \cite[Thm. 1.4]{BHP1} (translated from left- into right-action) can be stated as follows:
\begin{theorem}[Auto-correlation formula for model sets in groups]\label{BHP1Formula}
If $\Lambda = \Lambda(G, H, \Gamma, W)$ is a model set in $G$ and $\widehat{\nu}$ denotes the unique $G$-invariant probability measure on $\Omega_\Lambda^\times$, then the associated auto-correlation measure $\widehat{\eta}$ is the unique Radon measure on $G$ which satisfies
\[
\pushQED{\qed} \widehat{\eta}(f\ast f^*) = \|\mathcal P_\Gamma(f \otimes {\bf 1}_W)\|^2_{L^2(\Gamma\backslash(G \times H))} \quad \text{for all }f\in C_c(G).\qedhere \popQED
\]
\end{theorem}

Now let ${}_Kp_*\Lambda$ denote the associated weighted model set in $K\backslash G$. By Remark \ref{InvariantMeasureWeightedModelSet}
the unique $G$-invariant measure $\nu$ on $\Omega^\times_{{}_Kp_*\Lambda}$ is given by the push-forward of $\widehat{\nu}$ under the continuous factor map ${}_Kp_*:\Omega_\Lambda \to \Omega_{{}_Kp_*\Lambda}$. Its auto-correlation measure $\eta$ can thus be obtained from the auto-correlation measure $\widehat{\eta}$ of $\widehat{\nu}$ by the formula in Proposition \ref{AutocorNatural}. Using the explicit formula from Theorem \ref{BHP1Formula} for $\widehat{\eta}$ and denoting by ${}_Kp_K: G \to K \backslash G/K$ the canonical projection we obtain the following formula for $\eta$:
\begin{corollary}[Auto-correlation formula for weighted model sets]\label{Autocor1} If $\nu$ denotes the unique $G$-invariant probability measure on the weighted model set ${}_Kp_*\Lambda(G, H, \Gamma, W)$ and ${}_Kp_*\Lambda$ denotes the associated weighted model set in $K\backslash G$, , then the associated auto-correlation measure ${\eta}$ is the unique Radon measure on $K\backslash G/K$ which satisfies
\[
\pushQED{\qed} {\eta}(f\ast f^*) = \|\mathcal P_\Gamma({}_Kp_K^*f \otimes {\bf 1}_W)\|^2_{L^2(\Gamma\backslash(G \times H))}  \quad \text{for all }f\in C_c(K\backslash G/K).\qedhere \popQED
\]
\end{corollary}
For actual computations of the auto-correlation measure, the following formula is often the most useful one:
\begin{proposition}[Summation formula for the auto-correlation]\label{PropAutocor2} The auto-correlation measure $\eta$ from Corollary \ref{Autocor1} is the unique Radon measure on $K\backslash G/K$ with 
\[
\eta(f \ast f^*) = \sum_{(\gamma_1, \gamma_2) \in \Gamma} {}_Kp_K^*(f \ast f^*)(\gamma_1) \cdot ({\bf 1}_W \ast {\bf 1}_{W^{-1}})(\gamma_2)  \quad \text{for all }f\in C_c(K\backslash G/K).
\]
\end{proposition}
\begin{proof} Set $\phi := {}_Kp_K^*f$ and $r := {\bf 1}_W$. Since ${}_Kp_K$ is a $*$-homomorphism we have to show that
\begin{equation}
\| \mathcal P_\Gamma(\phi \otimes r))\|_{L^2(Y, m_Y)}^2 = \sum_{(\gamma_1,\gamma_2) \in \Gamma}
(\phi * \phi^*)(\gamma_1) \cdot (r * r^*)(\gamma_2).
\end{equation}
Set $F := \phi \otimes r$ and denote by $\widetilde{\mathcal P}_\Gamma F(g,h) :=  \mathcal P_\Gamma F(\Gamma(g,h))$ the lift of $\mathcal P_\Gamma F$ to a function on $G \times H$. We compute
\begin{eqnarray*}
 \sum_{(\gamma_1,\gamma_2) \in \Gamma}
(\phi * \phi^*)(\gamma_1) \cdot (r * r^*)(\gamma_2)  
&=&  \sum_{(\gamma_1,\gamma_2) \in \Gamma} \int_G \phi(g)\phi^*(g^{-1}\gamma_1) dm_G(g) \int_H r(h)r^*(h^{-1}\gamma_2) dm_H(h)\\
&=&\int_{G \times H} \phi\otimes r(g,h)  \sum_{(\gamma_1,\gamma_2) \in \Gamma}  \overline{\phi \otimes r((\gamma_1^{-1}g, \gamma_2^{-1}h))}\, dm_G\otimes m_H(g,h)\\
&=& \int_{G \times H} F(x) \overline{\widetilde{\mathcal P}_\Gamma F(x)} \, dm_G \otimes m_H(x).
\end{eqnarray*}
Now denote by $\mathcal F$ a fundamental domain for the left-action of $\Gamma$ on $G \times H$. We then have
\begin{eqnarray*}
 \sum_{(\gamma_1,\gamma_2) \in \Gamma}
(\phi * \phi^*)(\gamma_1) \cdot (r * r^*)(\gamma_2)  
&=& \sum_{\gamma \in \Gamma}\int_{\gamma \mathcal F} F(x) \overline{\widetilde{\mathcal P}_\Gamma F(x)} \, dm_G \otimes m_H(x)\\
&=&  \sum_{\gamma \in \Gamma}\int_{\mathcal F} F(\gamma x) \overline{\widetilde{\mathcal P}_\Gamma F(\gamma x)} \, dm_G \otimes m_H(x)\\
&=& \int_{\mathcal F} \widetilde{\mathcal P}_\Gamma F(\gamma x) \overline{\widetilde{\mathcal P}_\Gamma F(\gamma x)}dm_G \otimes m_H(x)\\
&=& \| \mathcal P_\Gamma(\phi \otimes r))\|_{L^2(Y, m_Y)}^2.
\end{eqnarray*}
The proposition follows.
\end{proof}

\section{Auto-correlation distributions for translation bounded measures in the hyperbolic plane}\label{SecSl2}

\subsection{The hyperbolic plane and the group ${\rm SL}_2(\R)$}

Throughout this section let $G := {\rm SL}_2(\R)$ and define elements of $G$ by the formulas
\[
k_\theta 
:=
\left(
\begin{matrix}
\cos 2\pi \theta& \sin 2\pi \theta \\ 
-\sin 2\pi \theta & \cos 2\pi \theta
\end{matrix}
\right),
\quad
a_t 
:=
\left(
\begin{matrix}
e^{t/2} & 0 \\ 
0 & e^{-t/2}
\end{matrix}
\right)
\qand 
n_u 
:=
\left(
\begin{matrix}
1 & u \\ 
0 & 1
\end{matrix}
\right).
\]
Then the maps $\R\to G$ given by $\theta \mapsto k_\theta$, $t \mapsto a_t$ and $u \mapsto n_u$ are group homomorphisms, and we denote their images by $K = {\rm SO}_2(\R)$, $A$ and $N$ respectively. Note that $K \cong S^1$ whereas $A \cong N \cong \R$.

Recall that the groups $G$ acts by fractional linear transformation on the upper half-plane, i.e.
\[
\begin{pmatrix} a&b\\c&d \end{pmatrix}.z := \frac{az+b}{cz+d} \quad \left(\begin{pmatrix} a&b\\c&d \end{pmatrix} \in G, \, z\in \C,\, {\rm Im}(z)>0\right),
\]
preserving the hyperbolic metric. This action is transitive and the stabilizer of $i$ is given by $K$. We thus have $G$-equivariant homeomorphisms
\[
K\backslash G \to G/K \to \bH^2, \quad Kg \mapsto g^{-1}K \mapsto g^{-1}.i.
\]
We may thus think of the hyperbolic plane as the proper homogeneous space $K\backslash G$ of $G$. 

Throughout this section, $\mu \in \mathcal T_b(G \curvearrowright K\backslash G)$ denotes a translation bounded measure on the hyperbolic plane. We assume that $\Omega^\times_{\mu}$ is uniformly locally bounded and that there exists a $G$-invariant probability measure $\nu$ on $\Omega^\times_{\mu}$. Under these assumptions the auto-correlation measure $\eta$ of $\nu$ can be defined and is a Radon measure on $K\backslash G/K$. In particular, $\mu$ could be a weighted model set in the hyperbolic plane and $\nu$ would then be the unique $G$-invariant probabiliy measure on its hull.

%so that $k_{\theta_1}k_{\theta_2} = k_{\theta_1 + \theta_2}$, $a_{t_1}a_{t_2} = a_{t_1+t_2}$ and $n_{u_1}n_{u_2} = n_{u_1+u_2}$. We thus obtain abelian subgroups
%\[
%K := \{k_\theta \mid \theta \in [0, 1)\} = {\rm SO}_2(\R), \quad A := \{a_t \mid t \in R\} \qand N := \{n_u \mid u \in \R\}
%\]
%of $G$ which satisfy $K \cong S^1$ and $A \cong N \cong \R$ and $\bH^2 = K\backslash G$.
We will need the following basic facts concerning ${\rm SL}_2(\R)$ (see \cite{Lang}). Multiplication induces a diffeomorphism $A \times N \times K \to G$ and thus every $g \in G$ can be written uniquely as
\begin{equation}
g = a_{t} n_{u} k_{\theta}.
\end{equation}
If $f \in L^1(G)$ and $F(t, u, \theta) := f(a_tn_uk_\theta)$, then we will normalize Haar measure on $G$ such that
\[
\int_G f(g) \, dm_G(g) = \int_{[0,1)} \int_{\bR} \int_{\bR} F(t,u,\theta) \, dt \, du \, d\theta.
\]
We thus obtain a $G$-invariant measure $dm_{A\backslash G}$ on $A\backslash G$ by setting
\[
\int_{A\backslash G} f(x)\, dm_{A\backslash G}(x) \quad = \quad \int_{[0.1)} \int_{\bR} f(An_uk_\theta)\, du \, d\theta \quad (f \in C_c(A\backslash G)).
\]
This measure is uniquely determined by the fact that for $f \in C_c(G)$ we have
\begin{equation}\label{mAG}
\int_G f\, dm_G \quad = \quad \int_{A\backslash G} \int_{\R} f(a_tx)\, dt \, dm_{A\backslash G}(Ax) \quad(f \in C_c(G)).
\end{equation}
The group $A$ normalizes $N$, and we have
\begin{equation}\label{AvsN}
a_{t}n_u a_{-t} = n_{e^{t}u} \quad \textrm{for all $t,u \in \bR$}.
\end{equation}
This implies in particular that for $f \in C_c(N)$ and $t\neq0$ we have
\begin{equation}\label{NIntegral}
\int_{-\infty}^\infty f(a_{-t}n_{-u}a_{t}n_{u})\, du = \int_{-\infty}^\infty f(n_{(1-e^{-t})u}) \, du = \frac{1}{|1-e^{-t}|}\int_{-\infty}^\infty f(n_u)\, du.
\end{equation}
%\begin{equation}\label{NIntegral}
%\int_{-\infty}^\infty f(a_tn_ua_{-t}n_{-u})\, du = \int_{-\infty}^\infty f(n_{(e^t-1)u}) \, du = \frac{1}{|e^t-1|}\int f(n_u)\, du.
%\end{equation}
The \emph{Weyl group} $W = N_K(A)/Z_K(A)$ acts on $A$ by conjugation. If we define 
\[
w := \begin{pmatrix} 0&1\\ -1&0 \end{pmatrix},
\]
then $N_K(A) = \langle w \rangle$, $Z_K(A) = \langle w^2 \rangle$ and $wa_tw^{-1} = a_{-t}$. Thus $W \cong \Z/2\Z$ and a function $f: A \to \C$ is $W$-invariant if and only if it is even, i.e.\ $f(a_t) = f(a_{-t})$. The diffeomorphism $c_w: G \to G$ given by $x \mapsto wxw^{-1}$ descends to a diffeomorphism 
\begin{equation}\label{cbarw}\overline{c}_w: A\backslash G \to A\backslash G, \quad Ax \mapsto Awxw^{-1}
\end{equation} and since Lebesgue measure on $\R$ is invariant under sign change it follows from \eqref{mAG} that $\overline{c}_w$ preserves the measure $m_{A\backslash G}$.

The inclusion $A \hookrightarrow G$ induces a homeomorphism $W\backslash A \to K\backslash G/K$, and if we define $\widehat{\iota}: G \to [1, \infty)$ by $g \mapsto \frac{1}{2} \trace(g^Tg)$, then $\widehat{\iota}$ is bi-$K$-invariant and induces homeomorphisms
\begin{equation}\label{Defiota}
\iota: K\backslash G/K \to [1, \infty) \qand \iota_A: W\backslash A \to [1, \infty).
\end{equation}
Explicitly, we have $\iota_A(\{a_{\pm t}\})= \cosh(t)$ and thus $\iota_A^{-1}(r) = \{a_{\pm \arcosh(r)}\}$.

\subsection{The Harish transform and its inverse}
We recall the definition and basic properties of the Harish transform on ${\rm SL}_2(\R)$. Our exposition follows \cite{Lang}, but we decided to spell out some formulas more explicitly.
\begin{definition} Let $f \in C_c(G)$. Then the \emph{Harish transform} $\bH f: A \to \C$ is given by
\[
\big(\bH f\big)(a_t) = e^{t/2} \int_{-\infty}^\infty f(a_t n_u) \, du.
\]
\end{definition}
\begin{lemma}[Properties of the Harish transform] The Harish transform has the following properties:
\begin{enumerate}[(i)]
\item For all $f \in C_c(G)$ we have $\bH(f^*) = (\bH f)^*$.
\item $\bH(C_c(G, K)) \subset C_c(A)^W$ and $\bH(C_c^\infty(G, K)) \subset C_c^\infty(A)^W$.
\item For all $f_1, f_2 \in C_c(G, K)$ we have $\bH(f_1 * f_2) = \bH f_1 * \bH f_2$.
\end{enumerate}
In particular, $\bH$ yields $*$-algebra homomorphisms
\[
\bH: C_c(G, K) \to C_c(A)^W \qand \bH: C_c^\infty(G, K) \to C_c^\infty(A)^W.
\]
\end{lemma}
\begin{proof} (i) For all $f \in C_c(G)$ and $t \in \R$ we have by \eqref{AvsN},
\begin{eqnarray*}
(\bH f^*)(a_t) 
&=& 
e^{t/2} \int_{\bR} \overline{f(n_{-u} a_{-t})} \, du \quad = \quad e^{t/2} \int_{\bR} \overline{f(a_{-t}(a_{t}n_{-u} a_{-t}))} \, du \\
&=&
e^{t/2} \int_{\bR} \overline{f(a_{-t} n_{-e^{t}u})} \, du \quad = \quad e^{-t/2} \int_{\bR} \overline{f(a_{-t} n_{u})} \, du \quad = \quad(\bH f)^*(a_t).
\end{eqnarray*}

(ii) It is clear from the formula that $\bH$ preserves smoothness. Now let $f \in C_c(G,K)$; given $t \in \R$ the function $\widetilde{\phi_t}$ on $G$ given by $\phi_t(x) = f(x^{-1}a_tx)$ is left-$A$-invariant and hence descends to a function $\phi_t$ on $A\backslash G$. Using bi-$K$-invariance of $f$ and \eqref{NIntegral} we obtain for every $t \neq 0$ the formula
\begin{eqnarray*}
\bH(f)(a_t) &=& e^{t/2}\int_{-\infty}^\infty  \int_0^1 f(k_{-\theta} a_t n_u k_\theta) \, d\theta\, du\\
&=&{e^{t/2}}\cdot{|1-e^{-t}|} \int_0^1\int_{-\infty}^\infty  f(k_{-\theta} a_t (a_{-t}n_{-u}a_tn_u)k_\theta) \, du\, d\theta\\
&=& |e^{t/2}-e^{-t/2}|  \int_0^1 \int_{-\infty}^\infty \phi_t(A n_u k_{\theta}) \,d u\, d\theta\\
&=&  |e^{t/2}-e^{-t/2}| \cdot m_{A\backslash G}(\phi_t),
\end{eqnarray*}
and this formula extends to $t=0$ by continuity. Since $a_{-t} = wa_tw^{-1}$ and $m_{A\backslash G}$ is invariant under the diffeomorphism $\overline{c}_w$ from \eqref{cbarw} we 
have $m_{A\backslash G}(\phi_t) = m_{A\backslash G}(\phi_{-t})$. Since also $|e^{t/2}-e^{-t/2}|$ is an even function, we deduce that $t\mapsto \bH(f)(a_t)$ is even.

(iii) For all $f_1, f_2 \in C_c(G, K)$ and $t \in \R$ we have
\begin{eqnarray*}
\bH(f_1 * f_2)(a_t) 
&=& 
e^{t/2} \int_{\bR} \big(f_1 * f_2)(a_t n_u) \, du \\
&=& e^{t/2} \int_{\bR} \Big( \int_{[0,1)}  \int_{\bR} \int_{\bR} f_1(a_{\tau} n_{v} k_{\xi}) f_2(k_{\xi}^{-1} n_{v}^{-1} a_{-\tau+t} n_u) \, d\tau \, dv \, d\xi \,  \Big) \, du \\
&=&
e^{t/2} \int_{\bR} \Big( \int_{\bR} \int_{\bR} f_1(a_{\tau} n_{v}) f_2(n_{v}^{-1} a_{-\tau+t} n_u) \, d\tau \, dv \Big) \, du \\
&=&
e^{t/2} \int_{\bR} \Big( \int_{\bR} \int_{\bR} f_1(a_{\tau} n_{v}) f_2(a_{-\tau + t} \big(a_{-(-\tau+t)} n_{v}^{-1} a_{-\tau+t} \big) n_u) \, d\tau \, dv\Big) \, du \\
\end{eqnarray*}
Since by \eqref{AvsN} for all $\tau, t, v \in \R$ we have $a_{-(-\tau+t)} n_{v}^{-1} a_{-\tau+t}  = n_{-e^{-(-\tau+t)}v}$, we deduce that
\begin{eqnarray*}
\bH(f_1 * f_2)(a_t) 
&=&
e^{t/2} \int_{\bR} \Big( \int_{\bR} \int_{\bR} f_1(a_{\tau} n_{v}) f_2(a_{-\tau + t} n_{u-e^{-(-\tau+t)}v}) \, d\tau \, dv\Big) \, du \\
&=&
e^{t/2} \int_{\bR} \Big( \int_{\bR} \int_{\bR} f_1(a_{\tau} n_{v}) f_2(a_{-\tau + t} n_{u-e^{-(-\tau+t)}v}) \, du \, dv \Big) \, d\tau \\
&=&
e^{t/2} \int_{\bR} \Big( \int_{\bR} \int_{\bR} f_1(a_{\tau} n_{v}) f_2(a_{-\tau + t} n_{u}) \, du \, dv \Big) \, d\tau\\
&=& 
\int_{\bR} \bH f_1(a_\tau) \, \bH f_2(a_{-\tau + t} \big) \, d\tau \\
&=& \big(\bH f_1  * \bH f_2 \big)(a_t).
\end{eqnarray*}
This finishes the proof.
\end{proof}
To see that the morphism $\bH: C_c^\infty(G, K) \to C_c^\infty(A)^W$ is actually an isomorphism and describe its inverse explicitly, it is convenient to relate the Harish transform to the more classical Abel transform. The homeomorphisms $\iota$ and $\iota_A$ from \eqref{Defiota} induce isomorphisms
\[
\iota^*: C_c(G, K) \to C_c([1, \infty)) \qand \iota_A^*: C_c(A)^W \to C_c([1, \infty)).
\]
Under these isomorphisms, both $C^\infty_c(G, K)$ and $C^\infty_c(A)^W$ are mapped onto the subspace $C_c^\infty([1, \infty)) \subset C_c([1, \infty))$ consisting of fuctions in $C_c^\infty(\R)$ with support contained in $[1,\infty)$.

\begin{remark}[Harish transform vs.\ Abel transform]\label{HvsA} The \emph{Abel transform} 
\[
\bA: C_c^\infty([1, \infty)) \to C_c^\infty([1, \infty)), \quad \bA\phi(r) := \int_{-\infty}^\infty \phi(r + u^2/2)\,du\]
is related to the Harish transform by the commutative diagram
\[\begin{xy}\xymatrix{
C_c^\infty(G, K) \ar[rr]^{\bH}&&C_c^\infty(A)^W\\
C_c^\infty([1, \infty)) \ar[u]^{\iota^*}\ar[rr]^{\bA}&&C_c^\infty([1, \infty))\ar[u]_{\iota_A^*}.
}\end{xy}\]
Indeed, since $\iota(a_t n_u) = \cosh(t) + \frac{1}{2}(ue^{t/2})^2$ we have for all $\phi \in C_c^\infty([1, \infty))$ and $t \in \R$
\[\bH(\iota^*\phi)(a_t) = e^{t/2}\int_{-\infty}^\infty \phi( \cosh(t) + (ue^{t/2})^2/2)\, du =  \int_{-\infty}^\infty \phi(\cosh(t) + u^2/2)\,du  = \iota_A^*\bA(\phi)(a_t).\]
\end{remark}
\begin{lemma}[Inversion of the Abel transform] The Abel transform is a linear isomorphism with inverse given by
\[
\bA^{-1}\psi(r) = \frac{-1}{2\pi} \int_{-\infty}^\infty \psi'(r+v^2/2) \, dv \quad (\psi \in C_c^\infty[1, \infty)).
\]
\end{lemma}
\begin{proof} For all $r\in [1, \infty)$ the substitutions $\left\{\begin{array}{ll}u = R\cos \theta\\ v = R\sin \theta \end{array} \right\}$ and $\xi = R^2/2$ yield
\begin{eqnarray*}
\bA(\bA^{-1}\psi)(r) &=& \frac{-1}{2\pi}\int_{-\infty}^\infty\int_{-\infty}^\infty \psi'(r+u^2/2+v^2/2)\,du \,dv \quad = \quad  \frac{-1}{2\pi}\int_0^{2\pi}\int_{0}^\infty \psi'(r+R^2/2)R\,dR \,d\theta\\
&=&-\int_0^\infty \psi'(r+\xi) d\xi \quad = \quad -\left.\psi(r+\xi)\right|_0^\infty \quad = \quad \psi(r) - \lim_{\xi \to \infty} \psi(r + \xi) \quad = \quad \psi(r),
\end{eqnarray*}
where the last equality holds since $\psi$ has compact support.
\end{proof}
\begin{corollary}[Invertibility of the smooth Harish transform]\label{HInverse} The Harish transform yields an isomorphism of $*$-algebras $\bH: C_c^\infty(G, K) \to C_c^\infty(A)^W$ with inverse given by ${\bH}^{-1} = \iota^* \circ \bA^{-1} \circ (\iota_A^{-1})^*$.\qed
%\[
%\bH^{-1}\psi(g) = 
%\]
\end{corollary}
Note that we can write out the formula for the inverse Harish transform explicitly as follows: If $\psi = (\iota_A^{-1})^*(\phi)$, then $\psi(r) = \phi(\{a_{\arcosh(r)}\})$, and hence
\[
\psi'(r) = \frac{\frac d{dr} \phi (a_{\pm \arcosh(r)})}{\sqrt{r^2-1}}.
\]
We conclude that
\[
{\bH}^{-1}(\phi)(g) = \left.\frac{-1}{2\pi} \int_{-\infty}^\infty \frac{\frac{d}{dr} \phi\left(a_{\arcosh(r+u^2/2)}\right)}{\sqrt{(r+u^2/2)^2-1}}\, du\right|_{r = \frac{1}{2}\trace(g^\top g)}.
\]

\subsection{The inverse Harish transform on distributions}
Given a smooth manifold $M$ we equip $C_c^\infty(M)$ with the topology of uniform convergence of all derivatives on compacta. We then denote by $\mathcal D(M)$ the space of distributions on $M$, i.e. the space of continuous linear functionals on $C_c^\infty(M)$.
We claim that the inverse Harish transform induces a map
\[
(\bH^{-1})^*: \mathcal D(G)^{K \times K} \to \mathcal D(A)^W, \quad (\bH^{-1})^*\xi(f) := \xi(\bH^{-1}f) \quad (\xi \in  \mathcal D(G)^{K \times K}, f \in C_c^\infty(A)^W).
\]
This amount to showing that the inverse Harish transform is continuous. In view of Corollary \ref{HInverse} this is equivalent to continuity of the inverse Abel transform $\bA^{-1}: C_c^\infty[1, \infty) \to C_c^\infty[1, \infty)$.
We will establish the following more precise estimate:
\begin{lemma}[Continuity of the inverse Abel transform]\label{HInverseCont} For every $n \geq 0$,
\[
\|\bA^{-1}\phi^{(n)}\|_\infty \leq \frac{2\sqrt 2}{\pi}\cdot \max\{\|\phi^{(n)}\|_\infty, \|\phi^{(n+1)}\|_\infty\} \quad (\phi \in C_c^\infty[1, \infty)).
\]
\end{lemma}
\begin{proof} For $\psi \in C_c^\infty[1, \infty)$ the substitution $v:= u^2/2$ yields
\[
\bA^{-1}\psi(r) = \frac{-1}{2\pi} \cdot 2 \int_{0}^\infty \psi'(r+u^2/2)\,du = \frac{-1}{2\pi} \cdot 2 \int_{0}^\infty \psi'(r+v)\,\frac{dv}{\sqrt{2v}} = \frac{-\sqrt 2}{2\pi}\int_0^\infty \frac{\psi'(r+v)}{\sqrt v}\,dv.
\]
Now, on the one hand,
\[
\left|\int_0^1 \frac{\psi'(r+v)}{\sqrt v}\,dv\right| \leq \|\psi'\|_\infty \cdot \int_0^1 \frac{dv}{\sqrt v} = 2 \cdot \|\psi'\|_\infty,
\]
and on the other hand by partial integration
\[
\int_1^\infty  \frac{\psi'(r+v)}{\sqrt v}\,dv = \left.\psi(r+v)\cdot\frac{1}{\sqrt v}\right|_{1}^\infty-\int_1^\infty \psi(r+v) \cdot \frac{-1}{2v^{3/2}}\,dv,
\]
hence 
\[
\left|\int_1^\infty  \frac{\psi'(r+v)}{\sqrt v}\,dv \right| \leq \|\psi\|_\infty + \frac{1}{2}\|\psi\|_\infty \int_1^\infty \frac{dv}{v^{3/2}} = 2\|\psi\|_\infty.
\]
We deduce that
\[
\|\bA^{-1}\psi\|_\infty \leq \frac{\sqrt 2}{2\pi}\cdot (2\|\psi'\|_\infty + 2\|\psi\|_\infty) \leq \frac{2\sqrt 2}{\pi}\max\{\|\psi\|_\infty, \|\psi'\|_\infty\},
\]
and applying this to $\psi := \phi^{(n)}$ yields the lemma.
\end{proof}

\subsection{The auto-correlation as a positive-definite distribution}
Recall that $\eta \in M^+(K\backslash G/K)$ denotes the auto-correlation measure of $\nu$. Given $f \in C^\infty_c(A)$, we define $f_W \in C_c^\infty(A)^W$ by
\[
f_W (a_t) := \frac{f(a_t)+f(a_{-t})}{2}.
\]
Then $\bH^{-1}(f_W) \in C^\infty_c(G, K)$ and thus ${}_K(\bH^{-1}(f_W))_K \in C_c(K\backslash G/K)$. We may thus apply the measure $\eta$ to this function and define
\[
\xi(f) := \eta({}_K(\bH^{-1}(f_W))_K).
\]
\begin{proposition} The map $\xi: C^\infty_c(A) \to \C$ has the following properties:
\begin{enumerate}[(i)]
\item $\xi$ is $W$-invariant and continuous, hence defines a distribution $\xi \in \mathcal D(A)^W$.
\item $\xi$ is evenly positive-definite in the sense that $\xi(f \ast f^*) \geq 0$ for all $f \in C^\infty_c(A)^W$.
\item $\xi$ determines the auto-correlation measure $\eta$ uniquely.
\end{enumerate}
\end{proposition}
\begin{proof} (i) $W$-invariance holds by construction, and continuity follows from Lemma \ref{HInverseCont}. (ii) follows from Remark \ref{AutocorPosDef} and the fact that $\bH^{-1}$ is a $*$-homomorphism, since for $f \in C_c(A)^W$ we have
\[
\xi(f \ast f^*) =  \eta({}_K(\bH^{-1}(f \ast f^*))_K) = \eta({}_K(\bH^{-1}f)_K \ast({}_K(\bH^{-1}f)_K)^*) \geq 0.
\]
(iii) follows from the fact that $\bH^{-1}$ is onto $C_c^\infty(G, K)$ and the latter is dense in $C_c(G, K)$.
\end{proof}
\begin{definition} The evenly positive-definite distribution $\xi \in \mathcal D(A)^W$ is called the \emph{auto-correlation distribution} of $\nu$.
\end{definition}
To summarize, we can consider the auto-correlation associated with a translation bounded measure in the hyperbolic plane either as a positive-definite Radon measure on $K\backslash G/K$ or equivalently as a evenly positive-definite distribution on $A \cong \R$.

\begin{remark}[Non-temperedness of the auto-correlation distribution]
Recall that a distribution $\xi \in \mathcal D(\R)$ is called \emph{tempered} if it can be extended to a continuous linear functional on the Schwartz space $\mathcal S(\R)$. In general there is no reason why the auto-correlation distribution would be tempered. The problem is that uniform model sets in the hyperbolic plane grow exponentially, and hence to define the periodization of a function over the hull of a weighted model set, one needs some form of exponential decay. While the inverse Harish transform does map the Schwartz space to the Schwartz space, the Harish transform of a Schwartz function will only be super-polynomially (rather than exponentially) decaying, and thus cannot be periodized.
\end{remark}

An important tool in the study of positive-definite distributions on $\R$ is given by the Fourier transform. For \emph{tempered} distributions there is Fourier inversion theorem which implies that a tempered distribution is uniquely determined by its \emph{real} Fourier transform. There is no such inversion theorem for non-tempered distributions, and we will see in the sequel to the present article that in order to fully determine the auto-correlation distribution we need to work with a suitable complex Fourier transform, which in the present case turns out to be related to the spherical Fourier transform associated with the Gelfand pair $(G, K)$.

\subsection{Beyond ${\rm SL}_2(\R)$} We have used the Harish transform to transform the auto-correlation measure into an evenly positive-definite distribution on $\R$. The underlying argument is not specific to the case of ${\rm SL}_2(\R)$ but applies in a much wider context:

If $G$ is any semisimple Lie group with finite center, then $G$ admits an Iwasawa decomposition $G = NAK$ with $K<G$ maximal compact and $A \cong \R^n$, and there is a Harish transform, which defines an isomorphism of $*$-algebras $\bH: C_c^\infty(K\backslash G/K) \to C_c^\infty(A)^W$, where the Weyl group $W = Z_K(A)/N_K(A)$ is a finite reflection group. One can verify that the inverse Harish transform is continuous, and hence descends to (non-tempered) distributions.

Consequently, if $\mu$ is a translation bounded measure in the Riemannian symmetric space $K\backslash G$ whose hull $\Omega_\mu^\times$ is uniformly locally bounded and admits a $G$-invariant probability measure $\nu$ with auto-correlation measure $\eta \in M^+(K\backslash G/K)$, then one obtains an evenly positive-definite distribution $\xi \in \mathcal D(A)^W$ by the formula
\[
\xi(f) := \eta({}_K(\bH^{-1}(f_W))_K),
\]
where
\[
f_W(a) := \frac{1}{|W|} \sum_{w \in W} f(w(a)). 
\]
As in the ${\rm SL}_2(\R)$-case, the auto-correlation measure is uniquely determined by this distribution, which justifies calling $\xi$ the \emph{auto-correlation distribution} of $\nu$. Since $A \cong \R^n$ we can view this distribution as an evenly positive-definite distribution on $\R^n$ with respect to a finite reflection group.

Actually, the story does not end here. There is a version of the Harish transform for semisimple algebraic groups over non-Archimedean local fields, called the \emph{Satake transform}, and thus auto-correlation distributions can also be defined in the non-Archimedean setting. 

We plan to return to both classes of examples in future work.

\section{Approximation of the auto-correlation for weighted regular model sets}\label{SecApprox}
Throughout this section let $\Lambda = \Lambda(G, H, K, \Gamma)$ be a regular model set (not necessarily uniform) in $G$ and let ${}_Kp_*\Lambda$ be the corresponding weighted model set in $K\backslash G$. We denote by $\nu$ the unique $G$-invariant probability measure on $\Omega^\times_{{}_Kp_*\Lambda}$ and by $\eta \in M^+(K\backslash G/K)$ its auto-correlation measure.

\subsection{Reminder of the amenable case}
Assume first that the group $G$ is amenable. We recall from Proposition \ref{AutocorNatural} that the auto-correlation measure $\eta$ satisfies
\[
\eta({}_Kf_K) = \widetilde{\eta}(f) \quad (f \in C_c(G, K)),
\]
where $\widetilde{\eta}$ is the auto-correlation measure of the unique $G$-invariant probability measure on $\Omega^\times_\Lambda$. In view of this formula, \cite[Cor. 5.4]{BHP1} yields the following formula for $\eta$:
\begin{corollary}[Approximation formula in the amenable case]\label{ApproxAmenable} Assume that $G$ is amenable. Then for every weakly admissible left-F\o lner sequence $(F_t)$ in $G$ the auto-correlation measure $\eta$ is given by the sampling limit
\[
\eta({}_Kf_K) = \lim_{t\to \infty} \frac{1}{m_G(F_t)}\sum_{x \in \Lambda \cap F_t} \sum_{y \in \Lambda} f(xy^{-1}) \quad (f \in C_c(G, K)).
\]
%In particular, $\eta$ is supported on the discrete set ${}_Kp_K(\Lambda\Lambda^{-1})$.
\end{corollary}
Let us recall the relevant definition of a weakly admissible sequence:
\begin{definition}\label{DefWeaklyAdmissible}
We say that a sequence $(F_t)_{t>0}$ of compact subsets of $G$ is 
\emph{weakly admissible} if each $(F_t)$ has positive Haar measure and there are continuous functions $\alpha, \beta : [0,1) \ra \bR_{+}$ 
with $\alpha(0) = 0$ and $\beta(0) = 0$ such that
\[
(F_t)_\delta \subset F_{t + \alpha(\delta)}
\qand
\sup_s \frac{m_G(F_{s + \delta})}{m_G(F_s)} = 1  + \beta(\delta),
\]
for all $t, \delta > 0$. %We shall refer to the pair $(\alpha, \beta)$ as the \emph{parameters} of 
$(F_t)$.
\end{definition}
This is a weakening of the notion of an \emph{admissible sequence} from \cite{GorodnikN-10}.

\subsection{An approximation theorem for non-amenable groups}
We would like to establish a version of Corollary \ref{ApproxAmenable} also for certain non-amenable groups. Since these groups do by definition not admit any F\o lner sequences, we need to find a different set of assumptions concerning $G$ and $(F_t)$. We are going to work in the context of Howe--Moore groups in the sense of the following definition:
\begin{definition} The lcsc group $G$ is called a \emph{Howe--Moore group} if it is non-compact and for every unitary $G$-representation $(V, \pi)$ with $V^G = \{0\}$ and all $u,v$ in $V$ we have
\[
\langle \pi(g)(u), v \rangle \to 0 \quad \text{as }g \to \infty.
\]
\end{definition}
Many non-amenable groups of interest have this property:
\begin{example}
\begin{enumerate}
\item ${\rm SL}_n(\R)$ and ${\rm SL}_n(\Q_p)$ are Howe--Moore groups.
\item More generally, if ${\bf G}$ is a semisimple algebraic group over a local field $k$, then ${\bf G}(k)$ is a Howe--Moore group  \cite{HM79}. In particular, every semisimple Lie group with finite center is a Howe--Moore group.
\item Finite (and also restricted infinite) products of Howe--Moore groups are Howe--Moore groups. In particular, products of real and $p$-adic semisimple algebraic groups and adelic semisimple algebraic groups are Howe--Moore groups.
\item If $G$ is the automorphism group of a regular tree, then $G$ is not a Howe--Moore group. However, $G$ contains a unique topologically simple subgroup $G_0$ of index $2$, which can be defined as the subgroup preserving any bi-partite coloring of the vertices (see \cite{Tits}), and this group is a Howe--Moore group (\cite{LubotzkyMozes}).
\item Similarly, automorphism groups of Bruhat--Tits buildings have a finite-index Howe--Moore subgroup.
\end{enumerate}
We will be interested in the following property of Howe--Moore groups:
\end{example}
\begin{remark}[Relative Howe--Moore property] Let $Y := \Gamma \backslash (G \times H)$ and denote by $L^2_0(Y)$ the orthogonal complement of the constants in $L^2(Y)$ so that the unitary representation $\pi: G \to \mathcal U(L^2_0(Y))$ has no invariant vectors. Moreover, let $(F_t)$ be a weakly-admissible sequence in $G$ and define $\beta_t := \frac{1}{m_G(F_t)} \cdot {\bf 1}_{F_t} \in L^1(G)$. We say that $G$ has the \emph{Howe--Moore property relative $(Y, K, (F_t))$} if 
\begin{equation}\label{RelHM}
\langle \pi(\beta_t) u, v \rangle \to 0 \quad \text{for all }u, v \in L^2_0(Y)^K.
\end{equation}
Note that if $G$ is a Howe--Moore group, then it has this property for any choice of $(Y, K, (F_t))$.
\end{remark}
In order to obtain an approximation theorem, we will impose an additional condition on the pair $(G, K)$. This condition is not strictly necessary to obtain an approximation theorem, but it is satisfied in all our examples of interest and simplifies the proof considerably.
\begin{definition} The pair $(G, K)$ is called a \emph{Gelfand pair} and $K\backslash G$ is called a \emph{commutative space} if the Hecke algebra $C_c(G, K)$ is commutative under convolution. 
\end{definition}
\begin{example} In developing the current theory we had the following examples of Gelfand pairs $(G, K)$ and commutative spaces $X = K\backslash G$ in mind (cf. \cite{Wolf}):
\begin{enumerate}
\item $G$ is abelian and $K = \{e\}$ so that $X = G$. This is the classical setting of quasi-crystals.
\item $G = \R^n \rtimes O(n)$, $K= O(n)$ so that $X= \R^n$ and $C_c(G, K)$ can be identified with radial functions on Euclidean space. This setting is implicitly studied in \cite{BFG}.
\item $G = {\rm SL}_2(\R)$ and $K = {\rm SO}_2(\R)$. In this case, $X$ can be identified with the hyperbolic plane and $C_c(G,K)$ can be identified with radial functions on the hyperbolic plane.
\item $G = {\rm SL}_2(\Q_p)$ and $K = {\rm SL}_2(\Z_p)$. In this case, $X$ can be identified with the vertex set of a $(p+1)$-regular tree, and $C_c(G, K)$ can be identified with radial functions on the vertices of the tree.
\item Generalizing (3), choose $G$ to be any semisimple Lie group with finite center and $K< G$ a maximal compact subgroup. In this case, $X$ is a Riemannian symmetric space. 
\item Generalizing (4),  choose $G = {\bf G}(k)$ and  $K = {\bf G}(o_k)$ where ${\bf G}$ is a semisimple algebraic group over a non-Archimedian local field $k$ and $o_k \subset k$ is its valuation ring. In this case $X$ can be identified with an orbit of special vertices in the Bruhat--Tits building of $G$. 
\item Another generalization of (4) is given as follows: Let $G$ be the automorphism group of a regular tree $T$ and let $K$ be the stabilizer of a vertex; in this case, $X$ is the vertex set of $T$.
\item Let $H$ be the $(2n+1)$-dimensional Heisenberg group. The group ${\rm U}(n)$ acts on $H$ by automorphisms fixing the center, and we can choose $K$ to be any subgroup of $U(n)$ containing a maximal torus and set $G := K \ltimes H$. Then $X$ can be identified with the Heisenberg group. If $K = {\rm U}(n)$, then $C_c(G, K)$ corresponds to radial functions on the Heisenberg group, and if $K$ is a maximal torus, then $C_c(G, K)$ corresponds to polyradial functions on the Heisenberg group.
\item One can consider finite products (and even restricted infinite products) of all of the pairs above. This includes in particular $S$-adic and adelic semisimple groups.
\end{enumerate}
Note that in all of these examples, either $G$ is amenable or $G$ has a finite index subgroup which is a Howe--Moore group.
\end{example}

We can now formulate a version of the approximation theorem for non-amenable groups:
\begin{theorem}[Approximation formula for Howe--Moore Gelfand pairs]\label{ApproximationThm} Assume that $G$ is a Howe--Moore group and that $(G, K)$ is a Gelfand pair. Then for every weakly admissible sequence $(F_t)$ of bi-$K$-invariant subsets in $G$ the auto-correlation measure $\eta$ is given by the sampling limit
\[
\eta({}_Kf_K) = \lim_{t\to \infty} \frac{1}{m_G(F_t)} \sum_{x \in \Lambda \cap F_t} \sum_{y \in \Lambda} f(xy^{-1}) \quad (f \in C_c(G, K)).
\]
\end{theorem}
As we will see, the theorem still holds if $G$ is not necessarily a Howe--Moore group, but has the Howe--Moore property relative to $(Y, K, (F_t))$.
\subsection{Proof of the approximation theorem}
Throughout this subsection we assume that $G, K, (F_t)$ satisfy the assumptions of the approximation theorem (Theorem \ref{ApproximationThm}). We denote $Y := \Gamma\backslash(G \times H)$ and $\beta_t := \frac{1}{m_G(F_t)}\cdot{\bf 1}_{F_t} \in L^1(G, K)$. By a standard argument, the (relative) Howe--Moore property implies the mean ergodic theorem, which we can formulate as follows:
\begin{lemma}[Mean ergodic theorem]\label{METConvenient} Let $\sigma$ be a $K$-invariant probability measure on $Y$ which is absolutely continuous with respect to $m_Y$ with square-integrable density. Then for every $\varphi \in C_b(Y)$,
\[
\check \beta_t \ast \sigma(\varphi) \to m_Y(\varphi).
\]
\end{lemma}
\begin{proof}
By assumption, $\sigma = u \, dm_Y$ with $u \in L^2(Y)^K$. For every $\varphi \in C_b(Y)$, 
\begin{eqnarray*}
(\check{\beta}_t * \sigma)(\varphi)
&=&
\int_G \check \beta_t(g)(g.\sigma)(\varphi) \, dm_G(g) = \int_G \beta_t(g)\sigma(g.\varphi) \, dm_G(g) \\
&=&
\int_G  \beta_t(g)\Big( \int_Y \varphi(g^{-1}.y) \, d\sigma(y) \Big) \, dm_G(g) = \int_{G}  \beta_t(g) \langle \pi_Y(g)\varphi,u \rangle \, dm_G(g)\\
&=& \langle \pi(\beta_t) \varphi, u \rangle
\end{eqnarray*}
Since $\beta_t$ is $K$-invariant, we have $\langle \pi(\beta_t) \varphi, u \rangle =  \langle \pi(\beta_t) \varphi^K, u \rangle$, where $\varphi^K$ denotes the projection of $\varphi$ onto $L^2(Y)^K$. If we write $\varphi^K = m_Y(\varphi) + \varphi_o$ and $u = 1 + u_o$ with $\varphi_o,u_o \in L^2_o(Y)$, then
\[
(\check{\beta}_t * \sigma)(\varphi) =  \langle \pi(\beta_t) \varphi^K, u \rangle =  m_Y(\varphi) +  \langle \pi(\beta_t) \varphi_o, u_o \rangle,\]
and since $G$ has the Howe-Moore property with respect to $(Y,K,(F_t))$, the last term tends to zero as $t \ra \infty$.
\end{proof}
The work of Gorodnik and Nevo \cite{GorodnikN-10, GorodnikN-122, GorodnikN-121, GorodnikN-14, GorodnikN-15} investigates in great generality under which conditions one can sharpen mean ergodic theorems for (possibly non-amenable) groups into pointwise statements. In many cases a pointwise statement can be derived without assuming commutativity of $C_c(G, K)$. However, since commutativity of $C_c(G, K)$ holds in essentially all examples of interest to us, we confine ourselves to this case, in which there is a particularly simple proof. Recall that the \emph{vague topology} on bounded measures on $Y$ is the weak-$*$-topology with respect to $C_0(Y)$, and that the subset of sub-probability measures is compact in this topology.

\begin{proposition}[Pointwise ergodic theorem]\label{PET}
In the situation of Theorem \ref{ApproximationThm} we have vague convergence
\[
\check \beta_t \ast \delta_y \to m_Y \quad \text{for every $y \in Y$.}
\]
\end{proposition}
\begin{proof} Since the space of sub-probability measures on $Y$ is vaguely compact, it suffices to show that every limit point $\nu  = \lim_{n \to \infty} \check \beta_{t_n} \ast \delta_y$ of $(\check \beta_t \ast \delta_y)$ coincides with $m_Y$. To show this, let $\rho = \rho_1 \otimes \rho_2 \in C_c(G \times H)$ be a probability density such that $\rho_1\in C_c(G, K)$ is bi-$K$-invariant. Since $\check \beta_t$ and $\rho_1$ commute, we then have 
\begin{eqnarray*}
\rho \ast \nu &=&  \lim_{n \to \infty}\, \rho \ast (\check \beta_{t_n} \ast \delta_y)\quad = \quad  \lim_{n \to \infty}\, ((\rho_1 \ast \check \beta_{t_n}) \otimes \rho_2) \ast \delta_y\\
&=&  \lim_{n \to \infty}\, ((\check \beta_{t_n} \ast \rho_1) \otimes \rho_2) \ast \delta_y \quad = \quad \lim_{n \to \infty} \check \beta_{t_n} \ast (\rho \ast \delta_y).
\end{eqnarray*}
Now let $\sigma := \rho \ast \delta_y$. We claim that $\sigma = \psi \cdot m_Y$ is absolutely continuous with respect to $m_Y$ with uniformly bounded density $\psi$. For the proof we may assume without loss of generality that $y = \Gamma$. If $F \in C_c(Y)$, then there exists $f \in C_c(G \times H)$ such that $F = \mathcal P_{\Gamma} f$, and thus
\begin{eqnarray*}
\sigma(F) &=& \quad  \int_G \int_H \sum_{(\gamma_1, \gamma_2) \in \Gamma}\rho(g,h) f((\gamma_1 g, \gamma_2 h)) \, dm_G(g)dm_H(h) \\
&=&  \int_G \int_H  \sum_{(\gamma_1, \gamma_2) \in \Gamma}\rho(\gamma_1^{-1}g,\gamma_2^{-1}h) f(g, h)\, dm_G(g)dm_H(h)\\
&=& \int_G \int_H \mathcal P_\Gamma \rho(\Gamma(g,h)) f(g, h)\, dm_G(g)dm_H(h)\\
&=& \int_Y \psi(y) F(y) dm_Y(y),
\end{eqnarray*}
where $\psi := \mathcal P_\Gamma \rho$. This proves that $\sigma = \psi \cdot m_Y$, and to see that this density is bounded, we observe that
\[
|\psi(\Gamma(g,h))| \leq \|\rho\|_\infty \cdot | \supp(\rho) \cap \Gamma(g,h)| %\leq C_K \|\psi\|_\infty, 
\]
Since the orbit of $\Gamma$ is uniformly locally finite and $\supp(\rho)$ is compact, the claim follows. We may thus apply Lemma \ref{METConvenient} to obtain
\[
\rho \ast \nu  =  \lim_{n \to \infty} \check \beta_{t_n} \ast (\rho \ast \delta_y) =  \lim_{n \to \infty} \check \beta_{t_n} \ast \sigma = m_Y.
\]
Now we can choose $\rho^{(n)}_1$ to be a convenient approximate identity in $C_c(G, K)$ and $\rho_2^{(n)}$ to be a convenient approximate identity in $C_c(H)$ and set $\rho^{(n)} := \rho^{(n)}_1 \otimes \rho^{(n)}_2$. Since $\nu$ is $K$-invariant by construction we then obtain
\[
\nu = \lim_{n \to \infty} \rho^{(n)} \ast \nu =  \lim_{n \to \infty}  m_Y = m_Y.
\]
This finishes the proof.
\end{proof}
We observe that by a standard approximation argument the convergence
\[
\check \beta_t \ast \delta_y(f) \to m_Y(f)
\]
holds not only for $f \in C_0(Y)$, but also for any compactly supported bounded function $f$ on $Y$, which is Riemann integrable with respect to $m_Y$. Theorem \ref{ApproximationThm} now follows from this pointwise statement and our previous work in \cite{BHP1}:
\begin{proof}[Proof of Theorem \ref{ApproximationThm}]
In \cite[Thm. 3.1]{BHP1} we introduced a parametrization map for the hull of a regular model set. Taking into account our change in convention from left- to right-actions, this map yields a $G$-equivariant surjection of the form $\beta: \Omega^\times_\Lambda \to Y$, where $Y = \Gamma\backslash (G \times H)$ as before, which induces an isomorphism
\[
\beta^*: L^2(Y) \to L^2(\Omega_\Lambda),
\]
where both $L^2$-spaces are with respect to the respective unique $G$-invariant measures. By \cite[Thm. 4.11 and Lemma 4.12]{BHP1} we have for all 
$f \in C_c(G)$ and for $m_G$-almost every $g \in G$,
\begin{equation}\label{EqA1}
\mathcal P_\Lambda f(\Lambda) = \mathcal P_\Gamma(f \otimes {\bf 1}_{W})(\beta(g.\Lambda)).
\end{equation}
In particular, if we abbreviate $y_0 := \beta(\Lambda)$, then for all $f\in C_c(G)$ and for $m_G$-almost every $s \in G$ we have
\[
 \mathcal P_\Lambda f(s^{-1}\Lambda) = \mathcal P_\Gamma(f \otimes {\bf 1}_{W})(\beta(s^{-1}\Lambda))  = \mathcal P_\Gamma (f \otimes {\bf 1}_{W})(s^{-1}y_0).
\]
Now let $f_1, f_2 \in C_c(G, K)$. We observe that the bounded measurable function
 \[h(y) :=  \mathcal P_\Gamma(f_1 \otimes{\bf 1}_W)(y) \overline{\mathcal P_\Gamma(f_2 \otimes{\bf 1}_W)(y)},\]
has compact support. Indeed, if $\mathcal F$ denotes a fundamental domain of $\Gamma$ in $G \times H$, then the supports of the functions $f_j \otimes{\bf 1}_W$ intersect only finitely many translates of $\mathcal F$, each of them in a compact set. Moreover, since $W$ is Jordan measurable, the function $h$ is Riemann integrable with respect to $m_Y$. Then \eqref{EqA1} yields
\[
\eta(f_1 \ast f_2^*) = \langle \mathcal P_\Lambda f_1, \mathcal P_\Lambda f_2 \rangle = \langle \mathcal P_\Gamma(f_1 \otimes {\bf 1}_W),  \mathcal P_\Gamma(f_2 \otimes {\bf 1}_W) \rangle = m_Y(h).
\]
Now if we set $\beta_t := \frac 1{m_G(F_t)}\cdot {\bf 1}_{F_t}$ the pointwise ergodic theorem (Proposition \ref{PET}) yields
\begin{eqnarray*}
\eta(f_1 \ast f_2^*) &=& m_Y(h) \quad = \quad \lim_{t \to \infty} \check \beta_t \ast \delta_{y_0}(h)\\
&=& \lim_{t \to \infty}  \frac{1}{m_G(F_t)} \; \int_{F_t}   \mathcal P_\Gamma(f_1 \otimes{\bf 1}_W)(s^{-1}.y_0) \overline{\mathcal P_\Gamma(f_2 \otimes{\bf 1}_W)(s^{-1}.y_0)}dm_G(s)\\
&=&\lim_{t \to \infty} \frac{1}{m_G(F_t)} \int_{F_t} \mathcal P f_1(s^{-1}\Lambda) \overline{\mathcal P f_2(s^{-1} \Lambda)} \, dm_G(s).
\end{eqnarray*}
In view of \cite[Theorem 5.3]{BHP1} this implies the theorem.
\end{proof}

\newpage

\appendix

\section{Convolution structures on double coset spaces}

\subsection{Convolution algebras and representations}\label{SecConvolutionAlgebras}
Recall that $M_b(G)$ denotes the Banach space of bounded complex Radon measures on the lcsc group $G$. 
\begin{remark}[$M_b(G)$ as a Banach-$*$-algebra]
$M_b(G)$ is a Banach-$*$-algebra with convolution product and involution respectively given by 
\[
\mu \ast \nu(f) = \int_G \int_G f(xy) d\mu(x)d\nu(y) \qand \int_G f d\mu^* = \int_G f(x^{-1}) d\overline{\mu}(x) \quad (f \in C_c(G)).
\]
For general $\mu, \nu$ in $M(G)$ the convolution $\mu \ast \nu$ is not well-defined, but if one of them has compact support, then $\mu \ast \nu$ can be defined by the same formula. Similarly, the involution $*$ can be extended by the same formula to all of $M(G)$.
\end{remark}
\begin{remark}[$L^1(G)$ as a Banach-$*$-algebra]
Our choice of Haar measure yields an embedding $L^1(G) \hookrightarrow M_b(G)$, $f \mapsto f \cdot m_G$. Under this embedding $L^1(G)$ is a Banach-$*$-subalgebra, and for $f, g \in L^1(G)$ we have
\[
f \ast g(x) = \int_G f(y)g(y^{-1}x) dm_G(y) \qand f^*(x) = \overline{f(x^{-1})} \quad (x, y \in G).
\]
By the same formula we obtain a convolution action of $L^1(G)$ on $L^p(G)$ for all $1 \leq p \leq \infty$ so that $\|f \ast g\|_p \leq \|f\|_1 \|g\|_p$. If $f \in L^1(G)$ and $g \in L^\infty(G)$, then $f \ast g$ is continuous.
\end{remark}
\begin{remark}[Extending representations]
If $\pi: G \to \mathcal U(V)$ is a unitary representation of $G$ on a Hilbert space $V$, then we obtain a $*$-representation of the algebra $L^1(G)$ by the formula
\[
\pi: L^1(G) \to \mathcal B(V), \quad \pi(f)(u) := \int_G f(g) \pi(g)u\, dm_G(g),
\]
where the integral can be interpreted in the weak sense.
\end{remark}
\begin{example}
If $\pi_L: G \to \mathcal U(L^2(G))$ denotes the left-regular representation, then for $f \in L^1(G)$ and $u \in L^2(G)$ we have
\[
\pi_L(f)(u)(x) = \int_G f(g) \pi_L(g) u(x)\, dm_G(g) = \int_G f(g) u(g^{-1}x)\, dm_G(g) = f \ast u(x),
\]
i.e. $\pi_L$ is the action by left-convolution.
\end{example}
\begin{example} If $\pi_R: G \to \mathcal U(L^2(G))$ denotes the right-regular representation, then for $f \in L^1(G)$ and $u \in L^2(G)$ we have
\begin{eqnarray*}
\pi_R(f)(u)(x) &=& \int_G f(g) \pi_R(g) u(x)\, dm_G(g) \quad = \quad \int_G f(g) u(xg)\, dm_G(g)\\
 &=& \int_G u(g)f(x^{-1}g)\, dm_G(g) \quad = \quad  \int_G u(g)\check f(g^{-1}x)\,dm_G(g) \quad = \quad u \ast \check f(x),
\end{eqnarray*}
i.e. $\pi_R(f)$ acts by right-convolution by $\check f$.
\end{example}

\subsection{Convolution algebras of bi-$K$-invariant functions}

\begin{remark}[Subalgebras of bi-$K$-invariant functions and measures]
The group $G$ acts on functions on $G$ by $L_gf(x) := f(g^{-1}x)$ and $R_gf(x) := f(xg)$, and dually on measures. We denote by $M_b(G,K) \subset M_b(G)$ and $L^1(G, K)\subset L^1(G)$ the Banach-$*$-subalgebras consisting of measures and function classes which are bi-$K$-invariant. The spaces $M(G, K)$, $C(G,K)$, $L^p(G, K)$ etc.\ are defined similarly, 
\end{remark}
\begin{definition} The dense $*$-subalgebra $C_c(G, K) := L^1(G, K) \cap C_c(G)$ is called the \emph{Hecke algebra} of the pair $(G, K)$. 
\end{definition}
Note that if ${}_Kp_K: G \to K\backslash G/K$ denotes the canonical projection, then pullback induces a bijection ${}_Kp_K^*: C_c(K\backslash G/K) \to C_c(G, K)$. By transport of structure, $C_c(K\backslash G/K)$ thus inherits the structure of a $*$-algebra from the Hecke algebra. We are going to describe this structure explicitly in Subsection \ref{DoubleCosetConvolution}, using a canonical measure on $K\backslash G/K$. To define this measure we need a variant of Weil's formula concerning integration on homogeneous spaces, which we will recall in the next subsection.
\subsection{A Weil formula for strongly proper actions}
Let $H$ be a unimodular lcsc group. We fix a Haar measure $m_H$ on $Y$.
\begin{definition} A lcsc $H$-space $Y$ is called a \emph{strongly proper $H$-space} if the action of $H$ on $Y$ is proper and the quotient $H\backslash Y$ is Hausdorff and paracompact with respect to the quotient topology. 
\end{definition}
The following version of Weil's formula can be found e.g.\ in \cite[Thm. 2.2]{Juestel}.
\begin{lemma}[Weil formula for strongly proper actions]\label{WeilFormula} Let $Y$ be a strongly proper $H$-space and let $\eta$ be an $H$-invariant Radon measure on $Y$. Then there exists a unique Radon measure $\underline{\eta}$ on $H\backslash Y$ such that for all $f \in C_c(Y)$,
\[
\pushQED{\qed} 
\int_{Y} f(y) d\eta(y) =\int_{H\backslash Y}  \left(\int_H f(hy) dm_H(h)\right) d\underline{\eta}(Hy).
\qedhere
\popQED
\]
\end{lemma}
The formula stated in \cite[Thm. 2.2]{Juestel} is actually more involved, but under our standing assumption that $H$ be unimodular it simplifies to the form above. We emphasize that $\underline{\eta}$ depends on our choice of Haar measure $m_H$ on $H$.

\subsection{Measures and convolutions on double coset spaces}\label{DoubleCosetConvolution}
Recall that $K$ denotes a compact subgroup of our lcsc group $G$ with Haar probability measure $m_K$. Since the left-action of $K$ on $G$ is strongly proper and preserves $m_G$, we can apply Lemma \ref{WeilFormula} with $Y := G$ and $H := K$. We deduce that there exists a unique Radon measure $m_{K\backslash G}$ on $K\backslash G$ such that for all $f \in C_c(G)$,
\[
\int_{G} f(g) \, dm_G(g) =\int_{K\backslash G}  \left(\int_K f(kg) dm_K(k)\right) dm_{K\backslash G}(Kg).
\]
Similarly, there exists a unique Radon measure $m_{G/K}$ on $G/K$ such that for all $f \in C_c(G)$,
\[
\int_{G} f(g) \, dm_G(g) =\int_{G/K}  \left(\int_K f(gk) dm_K(k)\right) dm_{G/K}(gK).
\]
Now observe that also the $(K\times K)$-action on $G$ given by $(k_1, k_2).g := k_1gk_2^{-1}$  is strongly proper and preserves $m_G$. Applying Lemma \ref{WeilFormula} with $Y := G$ and $H:= K \times K$ thus yields a unique Radon measure $m_{K\backslash G/K}$ on $K\backslash G/K$ such that for all $f \in C_c(G)$,
\[
\int_{G} f(g) dm_G(g) =\int_{K\backslash G/K}  \left(\int_K \int_K f(k_1gk_2) dm_K(k_1) dm_K(k_2)\right) dm_{K\backslash G/K}(KgK).
\]
\begin{definition} Let $f_1, f_2 \in C_c(K\backslash G/K)$. Then the \emph{convolution} of $f_1$ with $f_2$ is defined as the function $f_1 \ast f_2 \in C_c(K\backslash G/K)$ given by
\[
(f_1 \ast f_2)(KgK) = \int_{K\backslash G/K} f_1(KhK) \left(\int_K f_2(Kh^{-1}kgK) dm_K(k)\right) dm_{K\backslash G/K}(KhK).
\]
\end{definition}
There is no natural convolution structure on $C_c(K\backslash G)$ or $C_c(G/K)$, but we can define a convolution operation $C_c(K\backslash G) \times C_c(G/K) \to C_c(K\backslash G/K)$:
\begin{definition} Let $f_1 \in C_c(K\backslash G)$ and $f_2 \in C_c(G/K)$. Then the \emph{convolution} of $f_1$ with $f_2$ is defined as the function $f_1 \ast f_2 \in C_c(K\backslash G/K)$ given by
\[
(f_1 \ast f_2)(KgK)
= \int_{K\backslash G} f_1(Kh) \left(\int_K {f_2(h^{-1}k^{-1}gK)} dm_K(k) \right) dm_{K\backslash G}(Kh),
\]
\end{definition}
Also note that on $C_c(K\backslash G/K)$ we have a natural involution given by
\[
f^*(KgK) := \overline{f(Kg^{-1}K)} \quad (f \in C_c(K\backslash G/K)).
\]
Similarly, we have mutually inverse isomorphisms $C_c(K\backslash G) \to C_c(G/K)$ and  $C_c(G/K) \to C_c(K\backslash G) $ given by 
\[
f_1^*(gK) := \overline{f_1(Kg^{-1})} \qand f_2^*(Kg) := \overline{f_2(g^{-1}K)}  \quad (f_1 \in C_c(K\backslash G), f_2 \in C_c(G/K)).
\]
\begin{remark}[Relation to convolution in $C_c(G)$] Let us denote by ${}_Kp: G \to K\backslash G$, $p_K: G \to G/K$ and ${}_Kp_K: K\backslash G/K$ the canonical projections. We also denote by $C_c(G)^{L(K)}$ and $C_c(G)^{R(K)}$ the spaces of left-, respectively right-$K$-invariant functions in $C_c(G)$ so that $C_c(G, K) = C_c(G)^{L(K)} \cap C_c(G)^{R(K)}$. Since $K$ is compact we have bijections
\[
{}_Kp^*: C_c(K\backslash G) \to C_c(G)^{L(K)}, \; p_K^*: C_c(G/K) \to C_c(G)^{R(K)} \qand {}_Kp_K^*: C_c(K\backslash G/K) \to C_c(G,K).
\]
We recall that $C_c(G,K)$ is a $*$-subalgebra of the convolution algebra $C_c(G)$. Moreover, convolution induces a map $C_c(G)^{L(K)} \times C_c(G)^{R(K)} \to C_c(G, K)$, and the involution $*$ on $C_c(G)$ exchanges $C_c(G)^{L(K)}$ and $C_c(G)^{R(K)}$. Under the bijections above, these structures correspond to the convolution structures and involutions just defined: For $f_1, f_2, f \in C_c(K\backslash G/K)$ we have
\[
{}_Kp_K^*(f_1 \ast f_2) ={}_Kp_K^*f_1 \ast {}_Kp_K^*f_2 \qand ({}_Kp_K^*(f))^* =  {}_Kp_K^*(f^*),
\]
and for $f_1 \in C_c(K\backslash G)$ and $f_2 \in C_c(G/K)$ we have
\[
{}_Kp_K^*(f_1 \ast f_2) ={}_Kp^*f_1 \ast p_K^*f_2 \qand p_K^*(f_1^*) = ({}_Kp^*f_1)^* \qand {}_Kp^*(f_2^*) = (p_K^*f_2)^*.
\]
In particular, the $*$-algebras $C_c(K\backslash G/K)$ and $C_c(G, K)$ are isomorphic under ${}_Kp_K^*$.
\end{remark}
\subsection{Approximate identities for the Hecke algebra}
\begin{remark}[Canonical retractions] Define a map $M(G) \to M(G,K)$, $\mu \mapsto \mu^\sharp$ by \[
\int_G f(x) d\mu^\sharp(x) = \int_G\int_K\int_K f(k_1xk_2^{-1}) dm_K(k_1) dm_K(k_2)d\mu(x) \quad (f \in C_c(G)),
\]
For all $\mu \in M(G)$ and $f \in C_c(G, K)$ we then have $\mu(f) = \mu^\sharp(f)$, i.e. $\mu$ and $\mu^\sharp$ restrict to the same linear functional on $C_c(G, K)$. The map $\mu \mapsto \mu^\sharp$ restricts to a retraction of Banach-$*$-algebras $M_b(G) \to M_b(G,K)$ of norm $1$, and further to algebra retractions $L^1(G) \to L^1(G, K)$ and $C_c(G) \to C_c(G, K)$. Explicitly, for $f \in L^1(G)$ we have
\begin{equation}\label{fsharp}
f^\sharp(x) = \int_K \int_K f(k_1xk_2) dm_K(k_1) dm_K(k_2).
\end{equation}
The same formula also yields retractions $L^p(G) \to L^p(G, K)$ for $1 \leq p < \infty$, and we obtain continuous convolution actions of $L^1(G, K)$ on $L^p(G, K)$ such that $(f\ast g)^\sharp = f^\sharp \ast g^\sharp$. 
\end{remark}
\begin{remark}[Convenient approximate identities]\label{ConvenientApproximateIdentity} Let $U_n$ be a nested sequence of pre-compact identity neighbourhoods in $G$ with $\bigcap U_n = \{e\}$ (which exists since $G$ is second countable). We choose $\widetilde{\rho}_n \in C_c(G)$ with the following properties:
\[
\widetilde{\rho}_n \geq 0, \quad \widetilde{\rho}_n^* = \check{\widetilde{\rho}}_n = \widetilde{\rho}_n, \quad {\rm supp}(\widetilde{\rho}_n) \subset U_n \qand \int_G \widetilde{\rho}_n \, dm_G = 1.
\]
Then $\widetilde{\rho}_n \ast f$ and $f \ast \widetilde{\rho}_n$ converge to $f$ in the following sense (\cite[Prop.\ 2.44]{Folland}):
\begin{itemize}
\item If $1\leq p < \infty$, then convergence holds in $L^p$.
\item If $f \in C_c(G)$, then convergence holds uniformly, and hence if $f \in C(G)$, then convergence holds uniformly on compacta, and in particular pointwise. 
\end{itemize}
Now set $\rho_n := \widetilde{\rho}_n^\sharp$. Then we have convergence $\rho_n \ast f \to f^\sharp$ and $f \ast \rho_n \to f^\sharp$ in the same sense. In particular, $(\widetilde{\rho}_n)$ is a two-sided approximate identity in $C_c(G)$ and $L^1(G)$, and $(\rho_n)$ is a two-sided approximate identity in $C_c(G, K)$ and $L^1(G, K)$. We refer to these as \emph{convenient approximate identities}. Via the isomorphism $C_c(K\backslash G/K) \cong C_c(G, K)$ we also obtain a convenient approximate identity on $C_c(K\backslash G/K)$ with analogous properties.
\end{remark}
\begin{lemma}\label{ff*Dense}
\begin{enumerate}[(i)]
\item The subset $\{f \ast f^* \mid f \in C_c(G)\} \subset C_c(G)$ spans a dense subspace.
\item The subset $\{f \ast f^* \mid f \in C_c(G, K)\} \subset C_c(G, K)$ spans a dense subspace.
\item The subset $\{f \ast f^* \mid f \in C_c(K\backslash G/K)\} \subset C_c(K\backslash G/K)$ spans a dense subspace.
\end{enumerate}
\end{lemma}
\begin{proof} The span of $\{f \ast f^* \mid f \in C_c(G, K)\}$ contains all elements of the form $f \ast g^*$ with $f, g\in C_c(G, K)$ by polarization. Choosing a convenient approximate identity for $g$ then yields (ii) and hence (iii), and (i) follows from (ii) by choosing $K := \{e\}$.
\end{proof}

\bibliographystyle{abbrv}

\end{document}